\documentclass[oneside]{amsart}
\usepackage[T1]{fontenc}
\usepackage{mathtools,bbm,amsthm,amsmath,amsfonts,amssymb,color,epic}
\usepackage{multirow,comment}
\usepackage[export]{adjustbox}
\usepackage{tikz,pgfplots}
\usepackage[ruled]{algorithm2e}
\usepackage{enumitem}
\usepackage{subfig}
\usepackage{float}
\usepackage{xcolor}
\usepackage{graphicx}
\usepackage{hyperref}
\pgfplotsset{compat=newest}
\usetikzlibrary{arrows,decorations,backgrounds,positioning}
\usepackage{todonotes}
\theoremstyle{plain}
\usepackage{cite}
\newtheorem{theorem}{Theorem}[section]

\newtheorem{definition}[theorem]{Definition}
\newtheorem{remark}[theorem]{Remark}

\def\Letters{A,B,C,D,E,F,G,H,I,J,K,L,M,N,O,P,Q,R,S,T,U,V,W,X,Y,Z}
\makeatletter
\@for \@l:=\Letters \do{%
\expandafter\edef\csname\@l bb\endcsname{%
\noexpand\ensuremath{\noexpand\mathbb{\@l}}}%
\expandafter\edef\csname\@l bf\endcsname{%
{\noexpand\bf \@l}}%
\expandafter\edef\csname\@l cal\endcsname{%
\noexpand\ensuremath{\noexpand\mathcal{\@l}}}%
\expandafter\edef\csname\@l eu\endcsname{%
\noexpand\ensuremath{\noexpand\EuScript{\@l}}}%
\expandafter\edef\csname\@l frak\endcsname{%
\noexpand\ensuremath{\noexpand\mathfrak{\@l}}}%
\expandafter\edef\csname\@l rm\endcsname{%
{\noexpand\rm \@l}}%
\expandafter\edef\csname\@l scr\endcsname{%
\noexpand\ensuremath{\noexpand\mathscr{\@l}}}%
}

\newcommand{\bs}[1]{{\boldsymbol#1}}
\renewcommand{\d}{\operatorname{d\!}}
\newcommand{\Cor}{\operatorname{Cor}}
\newcommand{\Cov}{\operatorname{Cov}}
\newcommand{\DefField}{{\bs\chi}}

\newcommand{\isdef}{\mathrel{\mathrel{\mathop:}=}}

\newcommand{\Sref}{S_0}

\newcommand{\xref}{\hat{\bs x}}

\newcommand{\LBref}{\Delta_S}
\newcommand{\bigO}{\mathcal{O}}
\newcommand{\lowC}{{\bs L}}
\newcommand{\norm}[1]{\left\lVert#1\right\rVert}

\definecolor{navy}{RGB}{102,153,255}
\definecolor{tuerkis}{RGB}{51,153,204}

\title[IGA of diffusion problems on
random surfaces]
{Isogeometric analysis of diffusion problems on
random surfaces}

\author{Wei Huang}
\author{Michael Multerer}
\email{\{wei.huang,michael.multerer\}@usi.ch}

\address{
Euler Institute,
Universit\`a della Svizzera italiana,
Via alla Santa 1, 6962 Lugano
}

\begin{document}

\begin{abstract}
In this article, we discuss the numerical solution of diffusion equations
on random surfaces within the isogeometric framework.
We describe in detail, how diffusion problems on random surfaces
can be modelled and how quantities of interest may be derived. In particular,
we employ a low rank approximation algorithm for the high-dimensional
space-time correlation of the random solution based on an online singular
value decomposition, cp.\ \cite{brand2002incremental}. Extensive numerical
studies are performed to validate the approach. In particular, we consider
complex computational geometries originating from surface triangulations. The latter can be recast into the isogeometric context by transforming
them into quadrangulations using the procedure from \cite{pietroni2016tracing}
and a subsequent approximation by NURBS surfaces.
\end{abstract}
\maketitle
\section{Introduction}
Many problems in science and engineering 
can be modelled by partial differential equations (PDEs). 
Especially, PDEs on surfaces appear in a variety of applications,
such as computer graphics \cite{wu2008diffusion, tian2009segmentation,
crane2013geodesics}, chemical engineering \cite{bassett1978diffusion,
ehrlich1980surface, zhdanov2013elementary} or biology \cite{faraudo2002diffusion,
novak2007diffusion, elliott2012modelling} and there already exists a multitude
of computational methods for their numerical solution. In \cite{dziuk1988finite,
dziuk2013finite}, the surface finite element method has been developed aiming
at explicit surface representations by means of meshes. Alternatively, a level 
set based PDE solver is discussed in \cite{bertalmio2001variational,
bertalmio2000framework} and the closest point method, addressing closest point
representations of surfaces, is proposed in \cite{ruuth2008simple,
macdonald2008level,macdonald2010implicit}. A general solver for PDEs on point
clouds can be found in \cite{liang2013solving}. Due to the wide-spread use and
the effectiveness of Non-Uniform Rational Basis Spline (NURBS) based surfaces
in industry, an Isogeometric Analysis (IGA), see \cite{hughes2005isogeometric,
cottrell2009isogeometric}, based framework for the solution of PDEs on
surfaces was introduced in \cite{juttler2016numerical}. 

The striking advantage of IGA over the classical finite element method
is the use of an exact geometry representation. In particular,
the same basis is chosen to represent the geometry and the solution.
This way, in IGA, the apriori \(L^2\)-error and \(H^1\)-error satisfy the
optimal convergence rates \(\bigO(h^{p+1})\) and \(\bigO(h^p)\)
respectively for elliptic PDEs on surfaces, cf.\ \cite{bazilevs2006isogeometric,
da2011some,dede2015isogeometric}. Herein, \(h\) and \(p\) are the mesh size
and the polynomial order of the NURBS basis functions respectively. 
The \(L^2\)-error for the IGA of the parabolic problems is explored in \cite{zhu2017isogeometric}. In this article,
we shall adopt the IGA framework for solving diffusion PDEs on random surfaces.

As we have outlined, PDEs on deterministic surfaces are already well
explored.
In many applications, however, such as chemical engineering and biology,
the computational domain can not be measured exactly or may be subject
to manufacturing tolerances. Indeed, within this article, we will employ
an interpolation approach to recast arbitrary surface triangulations into
the isogeometric framework, which might itself be considered a source
of epistemic uncertainty. In such situations, the computational domain is
a major source of uncertainty that needs to be accounted for. Uncertainty 
quantification of  PDEs on random domains and surfaces has previously been
considered in \cite{xiu2006numerical,harbrecht2016analysis,CCX21,Dju21}
and in the context of IGA in \cite{beck2019iga,DHHM20}. To model random
surfaces, we follow the approach presented in \cite{harbrecht2016analysis}
based on the Karhunen-Lo\`eve expansion of random deformation fields and
exploit that the knowledge of this field at the surface is sufficient for
the modelling of the random surface itself, cp.\ \cite{DHHM20}.
Having the random deformation field at our disposal, quantities of interest
(QoI), such as the expectation and the correlation, can be defined with
respect to the reference configuration. As the computation of the
expectation is straightforward, we will particularly focus here on the
computation of the space-time correlation, whose size may easily exceed the
memory of a modern computer. Therefore, we devise a suitable low rank
approximation, based on the online singular value decomposition suggested in
\cite{brand2002incremental}. We remark that similar low rank approximation
methods have already been studied in the context of reduced basis
methods, see \cite{chen2017reduced,spannring2017weighted,
kazashi2020stability}. However, in contrast to these works,
our goal is not to devise a surrogate model for the underlying PDE, but
rather to find an efficient means to represent the correlation.

The rest of this article is structured as follows. In
Section~\ref{sec:ProbForm}, we introduce the diffusion problem and the
random model under consideration. In Section~\ref{sec:IGA}, we discuss
the isogeometric representation of random domains and briefly recall the algorithm from
\cite{pietroni2016tracing}
for transferring surface triangulations into quadrangulations.
Section~\ref{sec:Varf} then discusses the variational formulation. Especially,
thanks to the isogeometric setting we can compute realizations of the solution
directly in the spatial configuration. Section~\ref{sec:Varf} is devoted to
the computation of quantities of interest, particularly to the low rank
approximation of the solution's correlation. In 
Section~\ref{sec:Numerics}, extensive numerical studies are presented and
concluding remarks are stated in Section~\ref{sec:Conclusion}.

\section{Problem formulation}\label{sec:ProbForm}
Let \((\Omega,\Fcal,\Pbb)\) denote a complete probability space,
where \(\Omega\) is the sample space, \(\Fcal\) denotes the
\(\sigma\)-field of events, and \(\Pbb\) is a probability measure. As model problem,
we consider the diffusion equation on a closed
random surface, i.e.,
\begin{align}\label{eq:heatuq}
\begin{cases}
\partial_t u(\omega, {\bs x},t)
- \LBref u(\omega, {\bs x},t) = f(\omega,{\bs x}), &\quad {\bs x}
\in S(\omega)\\
u(\omega,{\bs x},0) = u_0(\omega, {\bs x}), &\quad {\bs x} \in S(\omega)\\
\end{cases}
\quad\text{{for~\(\mathbb{P}\)-a.e.}\ }\omega\in\Omega.
\end{align}
Herein, \(\LBref\) denotes the Laplace-Beltrami operator,
 \(f\) is a source on the surface, and \(u_0\) denotes
 the initial condition. As the surface is
 closed, we do not need to impose any boundary conditions.
 Uniqueness of the solution is then obtained by considering
 mean zero functions. For the treatment of non-trivial
 boundary conditions in the random domain case, we refer to
 \cite{GP18}.
 
For the modelling of the random surface, we assume the
existence of a Lipschitz continuous
reference surface \(\Sref\) and a random deformation field 
\[
\DefField\colon\Omega \times \Sref \to \mathbb{R}^3
\]
such that there exists a constant \(C_{\operatorname{uni}} > 0\) with 
\begin{equation}\label{eq:unifVfield}
\|\DefField(\omega)\|_{C^1(\Sref;\mathbb{R}^3)},
\|\DefField^{-1}(\omega)\|_{%
C^1(S(\omega);\mathbb{R}^3)}
\leq C_{\operatorname{uni}}
\quad\text{{for~\(\mathbb{P}\)-a.e.}\ }\omega\in\Omega
\end{equation}
and 
\[
S(\omega) = \DefField(\omega, \Sref) \quad
\text{{for~\(\mathbb{P}\)-a.e.}\ }\omega\in\Omega.
\]

\begin{figure}[htb]
\begin{center}
\begin{tikzpicture}
\node (bunny_ref) at (-3.0,0) {\includegraphics[scale=0.15, clip=true,trim = 800 100 800 150]{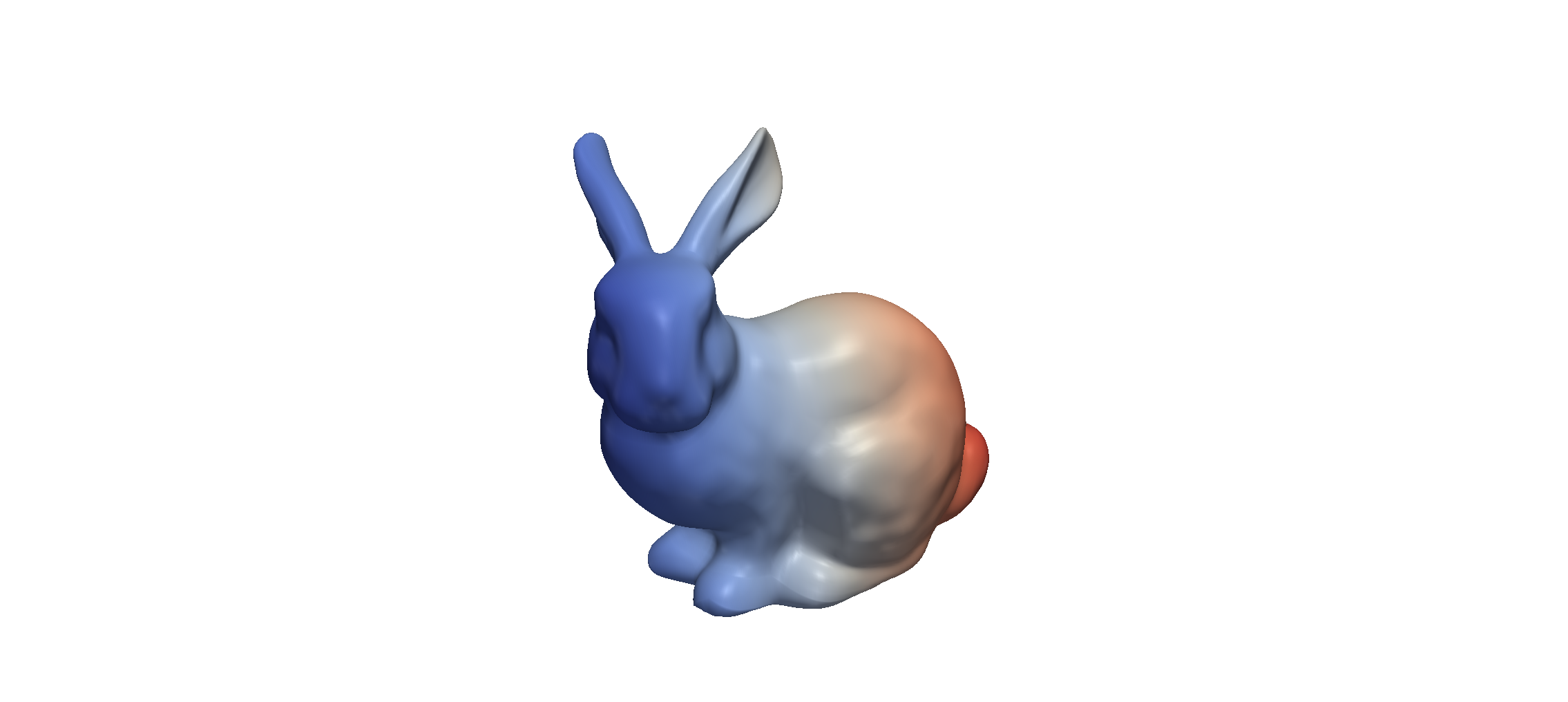}};
\node (bunny_def) at (3.0,0) {\includegraphics[scale=0.15, clip = true, trim = 800 0 750 150]{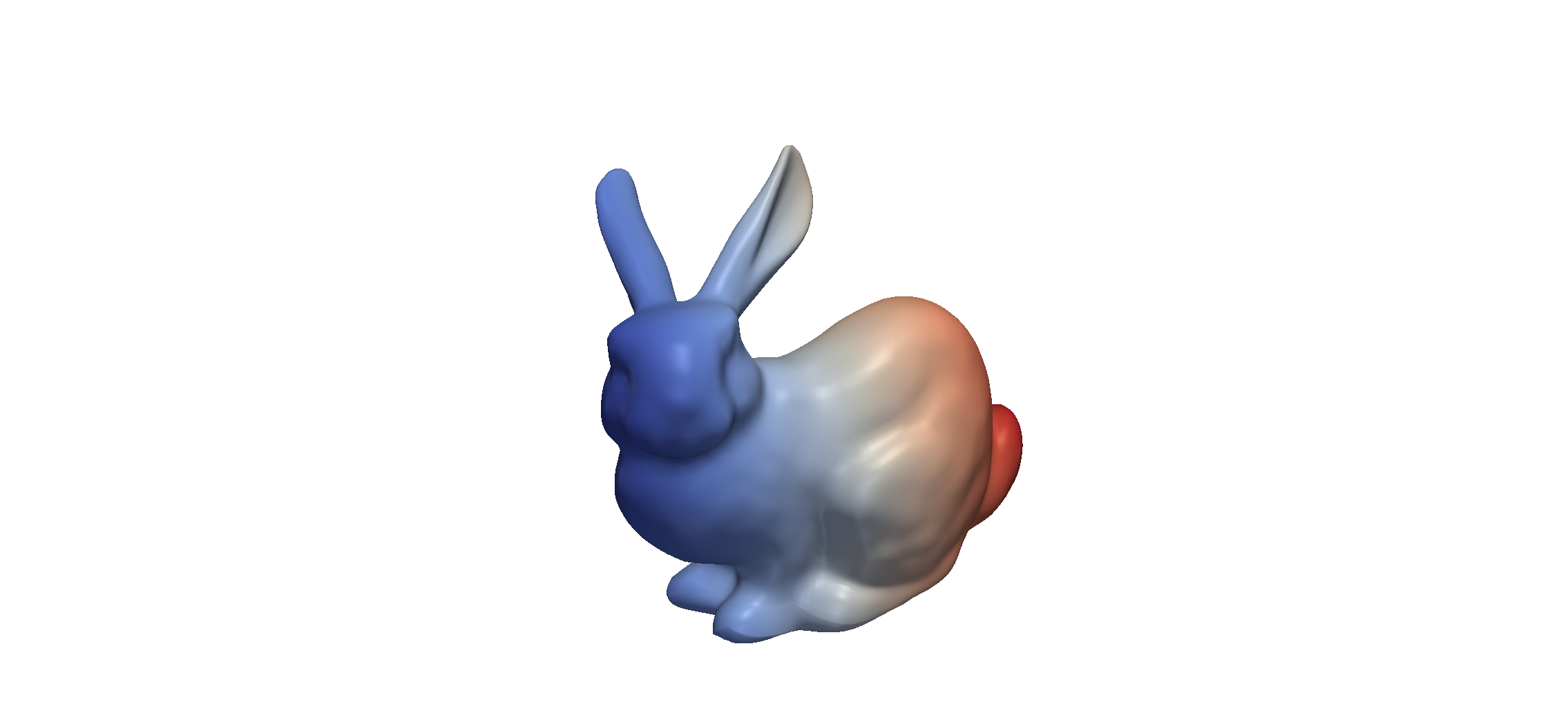}};
\node (a) [] at (-3.5,-1) {};
\node (b) [] at (3,-1) {};
\draw [ultra thick,black,->] (a) to[bend left] node[above]
(TextNode) {$\DefField(\omega, \hat{{\bs x}})$} (b);
\end{tikzpicture}
\caption{\label{fig:deformation}Realization of a randomly deformed surface.}
\end{center}
\end{figure}
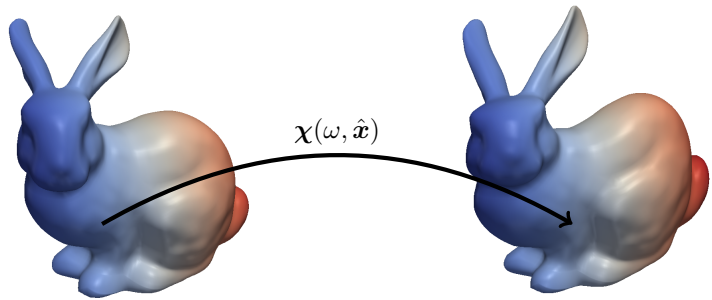

We refer to Figure~\ref{fig:deformation} for a realization of such
a random deformation field.
Note that it would in principle also be possible to consider
random space-time deformation fields, as they have recently
been introduced in \cite{GKMP20}, in a similar fashion.

\begin{remark}
For a sufficiently smooth surface, e.g.\ Lipschitz continuous,
the uniformity condition \eqref{eq:unifVfield} guarantees that
the Sobolev spaces \(H^1\big(S(\omega)\big)\) and \(H^1\big(S_0\big)\)
are equivalent independently of the particular realization of the
random parameter, see for example \cite{harbrecht2016analysis}.
A more general framework is considered in \cite{Dju21}, where
the bounds on the deformation field hold only pathwise. The latter
particularly covers the situation of log-normal deformation fields.
\end{remark}

For our model \eqref{eq:heatuq}, we make the assumption
that the initial condition
and the heat source on the random
surface are given in material coordinates
and not subject to uncertainty
themselves, i.e.,
\[
u_0(\omega, {\bs x})
\isdef\hat{u}_0\big(\DefField^{-1}(\omega,{\bs x})\big),\quad
f(\omega, {\bs x})
\isdef\hat{f}\big(\DefField^{-1}(\omega,{\bs x})\big)
\]
for two deterministic functions \(\hat{f},\hat{u}_0\colon S_0\to\Rbb\).
We remark that it is also possible to consider initial conditions and heat sources described
in spatial coordinates. In this case, they need to be 
defined on the hold-all domain \(\bigcup_{\omega\in\Omega}
\DefField(\omega,S_0)\). Furthermore, one may to consider the
situation where the data are subject to uncertainty themselves,
which typically results in a much higher dimensionality of the problem,
see e.g.\ \cite{Mul19}.
Finally, for a given solution \(u\) on \(S(\omega)\),
we denote the pulled back solution by
\[
\hat{u}(\omega, \xref,t)
\isdef u\big(\DefField(\omega, \xref),t\big).
\]

Based on the pulled back solution, we compute
QoI's, such as the expected solution and its space-time correlation. 
More precisely, for any fixed time points
$t$ and $t'$, we define these quantities according to
\begin{equation}
\mathbb{E}[\hat{u}(\omega,\xref,t)] = \int_\Omega \hat{u}(\omega, \xref,t) \d\mathbb{P}(\omega) \in H^1(\Sref)
\end{equation}
and 
\begin{equation}
\Cor[\hat{u}]\big((\xref,t),(\xref',t')\big) 
= \int_\Omega \hat{u}(\omega,\xref,t)
\hat{u}(\omega,\xref',t') \d\mathbb{P}(\omega)
\in H^1(\Sref)\otimes H^1(\Sref),
\end{equation}
respectively.

\section{Isogeometric representation of random domains}\label{sec:IGA}
\subsection{Representation of random deformation fields}
In this section we briefly recall the ideas presented in \cite{DHHM20}
and adapt them to our settings.
As has been argued there, it is sufficient to compute
the Karhunen-Lo\`eve expansion of the random deformation
field exclusively with respect to the random surface and no
volume discretization is required at all.
Hence, given the expected deformation field
\[
\Ebb[\DefField]\colon S_0\to\Rbb^d
\]
and its matrix valued covariance function
\[
\Cov[\DefField]\colon S_0\times S_0\to\Rbb^{d\times d},
\]
we can compute the surface Karhunen-Lo\`eve expansion
\begin{align}\label{eq:boundarykl}
\DefField(\omega,\xref)=\Ebb[\DefField](\xref)+
\sum_{k=1}^\infty\sqrt{\lambda_{k}}\DefField_k(\xref)
Y_{k}(\omega),\quad\xref\in\Sref.
\end{align}
Herein, the tuples \(\{(\lambda_{k}, \DefField_{k})\}_k\)
are the eigenpairs of the covariance operator
\begin{equation}\label{eq:surfacecovarianceoperator}
\begin{aligned}
&\mathcal{C}\colon [L^2(\Sref)]^3\to [L^2(\Sref)]^3,\\
&(\mathcal{C}{\bs v})(\xref)\isdef\int_{\Sref}
\Cov[\DefField](\xref,\xref')
{\bs v}(\xref')\d\sigma_{\xref'}
\end{aligned}
\end{equation}
and, for \(\lambda_k\neq 0\),
the centred and uncorrelated random variables
\(\{Y_k\}_k\) are given according to
\[
Y_{k}(\omega)\isdef\frac{1}{\sqrt{\lambda_{k}}}\int_{\Sref}
\big(\DefField({\omega},\xref)
-\Ebb[\DefField](\xref)\big)^\intercal\DefField_{k}(\xref)
\d{\sigma_{\xref}}.
\]
The existence of the Karhunen-Lo\`eve
expansion is guaranteed by the uniformity condition
\eqref{eq:unifVfield}, which ensures
\(\chi\in L^2\big(\Omega;[L^2(\Sref)]^3\big)\)
and \(\Cov[{\bs\chi}]\in[L^2(\Sref)]^3\otimes[L^2(\Sref)]^3\)
and thus gives rise to the Hilbert-Schmidt operator
\(\Ccal\) from \eqref{eq:surfacecovarianceoperator}, see also \cite{harbrecht2016analysis}. For the numerical realization
of the Karhunen-Lo\`eve expansion in the random domain
case, we refer to \cite{Mul19,HMvR21}.

In practice, however, the random variables \(\{Y_k\}_k\)
are not known explicitly and need to be estimated.
We make the common model assumption that the random variables
\(\{Y_{k}\}_k\) are independent and uniformly distributed with
\(\mathcal{U}(-1,1)\) for all \(k\). 
For numerical computations, the Karhunen-Lo\`eve expansion
has to be truncated after \(m\in\Nbb\) terms, where \(m\)
has to be chosen to meet the overall accuracy.
Then, by
identifying each random variable \(Y_k\) by its image
\(y_k\in[-1,1]\), we
arrive at the parametric deformation field
\begin{equation}\label{eq:parametricVectorFieldBoundary}
\DefField({\bs y},\xref)=\Ebb[\DefField](\xref)+
\sum_{k=1}^m\sqrt{\lambda_{k}}\DefField_{k}(\xref)y_k,
\quad{\bs y}\in\Gamma\isdef[-1,1]^{m}.
\end{equation}
The parametric deformation field gives rise to
the parametric surfaces
\begin{align}\label{eq:randomdomain}
  S({\bs y}) = \big\{\DefField({\bs y},\xref):
\xref\in \Sref\big\}.
\end{align}

In the next subsection, we define corresponding
NURBS representations and the associated ansatz spaces.
\subsection{Random NURBS surfaces}
We start by recalling the basic notions of isogeometric analysis,
restricting ourselves to function spaces constructed via locally
quasi-uniform $p$-open knot vectors as in \cite{BDK+2020}.

\begin{definition}\label{def::splines}
Let $p$ and $k$ with $0\leq p< k$. 
A \emph{locally quasi uniform $p$-open knot vector} is a tuple
\begin{align*}
\Xi = \big[{\xi_0 = \cdots =\xi_{p}}\leq \cdots \leq{\xi_{k}
=\cdots =\xi_{k+p}}\big]\in[0,1]^{k+p+1}
\end{align*}
with $\xi_0 = 0$ and $\xi_{k+p}=1$
such that there exists a constant $\theta\geq 1$ with
$\theta^{-1}\leq h_j\cdot h_{j+1}^{-1} \leq \theta$
for all $p\leq j < k$, 
where  $h_j\isdef  \xi_{j+1}-\xi_{j}$.
The B-spline basis $ \lbrace b_j^p \rbrace_{0\leq j< k}$ is then recursively
defined according to
\begin{align*}
b_j^p(x) & =\begin{cases}
\mathbbm{1}_{[\xi_j,\xi_{j+1})}&\text{ if }p=0,\\[8pt]
\frac{x-\xi_j}{\xi_{j+p}-\xi_j}b_j^{p-1}(x)
+\frac{\xi_{j+p+1}-x}{\xi_{j+p+1}-\xi_{j+1}}b_{j+1}^{p-1}(x)
& \text{ else,}
\end{cases}
\end{align*}
where $\mathbbm{1}_A$ refers to the characteristic function of the set $A$.
The corresponding spline space is defined as
$\Scal^p(\Xi)\isdef\operatorname{span}(\lbrace b_j^p\rbrace_{j <k}).$
\end{definition}

Spline spaces in two spatial dimensions are obtained by a
tensor product construction.
More precisely, we define the spaces 
\begin{equation}\label{eq:AnsatzSpace}
\Scal^{\bs p}\isdef \Scal^{p}(\Xi)\otimes
\Scal^{p}(\Xi).
\end{equation}
With regard to the knot vector $\Xi$,
sets of the form $[\xi_{i},\xi_{i+1}]\times[\xi_{j},\xi_{j+1}]$
are referred to as \emph{elements}. The corresponding mesh size
will be denoted by $h$.
For further details, we refer to \cite{piegl1996nurbs} and the references therein.

In what follows, we shall adopt the usual isogeometric setting for the
representation of the surface $\Sref$, i.e.,
we assume that $\Sref$ can be decomposed 
into several smooth \emph{patches}
\(
\Sref = \bigcup_{i=1}^M\Sref^{(i)},
\)
where the intersection $\Sref^{(i)}\cap\Sref^{(i')}$ 
consists at most of a common vertex or a common edge for 
\(i\neq i^\prime\).
Each patch $\Sref^{(i)}$ is represented by
an invertible NURBS mapping
\begin{equation}\label{eq:parametrization}
{\bs s}_i\colon[0,1]^2\to \Sref^{(i)}\quad\text{ with }
\quad \Sref^{(i)} = {\bs s}_i([0,1]^2)
\quad\text{ for } i = 1,2,\ldots,M.
\end{equation}
Herein, the functions \({\bs s}_i\) are of the form
\begin{align*}
{\bs s}_i(x_1,x_2)\isdef \sum_{0=i_1}^{k_1}
\sum_{0=i_2}^{k_2}\frac{{\bs c}_{i_1,i_2}
b_{i_1}^{p}(x_1) b_{i_2}^{p}(x_2) w_{i_1,i_2}}{
\sum_{j_1=0}^{k_1-1}\sum_{j_2=0}^{k_2-1} b_{j_1}^{p}(x_1)
b_{j_2}^{p}(x_2) w_{j_1,j_2}}
\end{align*}
for control points ${\bs c}_{i_1,i_2}\in \Rbb^3$ and weights
$w_{i_1,i_2}>0$. As a consequence, the resulting
patches \(S_0^{(i)}\) are at most of class \(C^{p-1}\). 
In order to obtain a consistent surface representation,
we moreover assume that parametrizations sharing a common edge
coincide at this edge except for orientation.

The subdivision of the surface \(S_0\) into
patches, directly induces a corresponding subdivision of the
random surface according to
\(S({\bs y})=\bigcup_{i=1}^M S^{(i)}(\bs y)\)
with \(S^{(i)}(\bs y)\isdef\big(\DefField({\bs y})\circ{\bs s}_i\big)([0,1]^2)\)
. By representing the mappings
\(\DefField({\bs y})\circ{\bs s}_i\) using NURBS, it can
be achieved that \(S^{(i)}(\bs y)\)
is again a NURBS surface, see \cite{DHHM20}. As
we typically have to compute many realizations of the parametric
surface to obtain acceptable approximations of the
QoI's, it is in practice advisable to reinterpolate each patch
and to compute the Karhunen-Lo\`eve expansion only with respect
to a set of interpolation points on \([0,1]^2\), see \cite{HMvR21}.

As usual, we define the ansatz space on \(\Sref\)
by patchwise lifting the space \(\Scal^{\bs p}\) to \(S_0\), i.e.,
\(
V_h\isdef
\big\{
v\in C(\Sref): \big(v|_{\Sref^{(i)}}\circ{\bs s}_i\big)\in\Scal^{\bs p}
\big\}.
\)
Analogously, we introduce the parametric spaces
\(
V_h({\bs y})\isdef
\big\{
v\in C\big(S({\bs y})\big): \big(v|_{\Sref^{(i)}}
\circ \DefField
({\bs y})\circ{\bs s}_i \big)\in\Scal^{\bs p}
({\bs\Xi})
\big\}.
\)
Given a B-spline
basis \(\widehat{\bs\Phi}\isdef[\hat{\varphi}_1,\hat{\varphi}_2,\ldots,
\hat{\varphi}_{N_{{\bs p},{\bs\Xi}}}]\) of \(V_h\), we denote the lifted basis in 
\(V_h({\bs y})\) by
\[\varphi_i\isdef(\hat{\varphi}_i\circ\DefField^{-1})\quad\text{for }
i=1,2,\ldots,N_{{\bs p},{\bs\Xi}}.
\]
\subsection{Partitioning surfaces into quadrilateral patches}
In practice, the surface representation \eqref{eq:parametrization}
is often not directly available for more complex geometries.
However, the partitioning of surface meshes into conforming
quadrilateral patches is well solved, see for example
\cite{campen2017partitioning, born2021layout, lyon2021quad,
lyon2021simpler,pietroni2016tracing} and the references therein.
The procedure suggested there is based on determining a suitable set
of nodes called \emph{singularities} and connecting these by a set
of arcs called \emph{separatrices}, such that the resulting partition
consists of multiple valid conforming quadrilateral patches. 

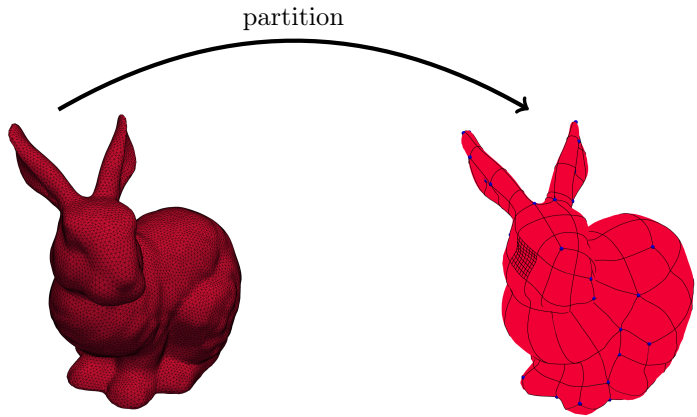
\begin{figure}[htb]
\begin{center}
\begin{tikzpicture}
\node (bunny_origin) at (-3.0,0) {\includegraphics[
width=0.3\textwidth, clip=true,trim = 100 0 100 0]{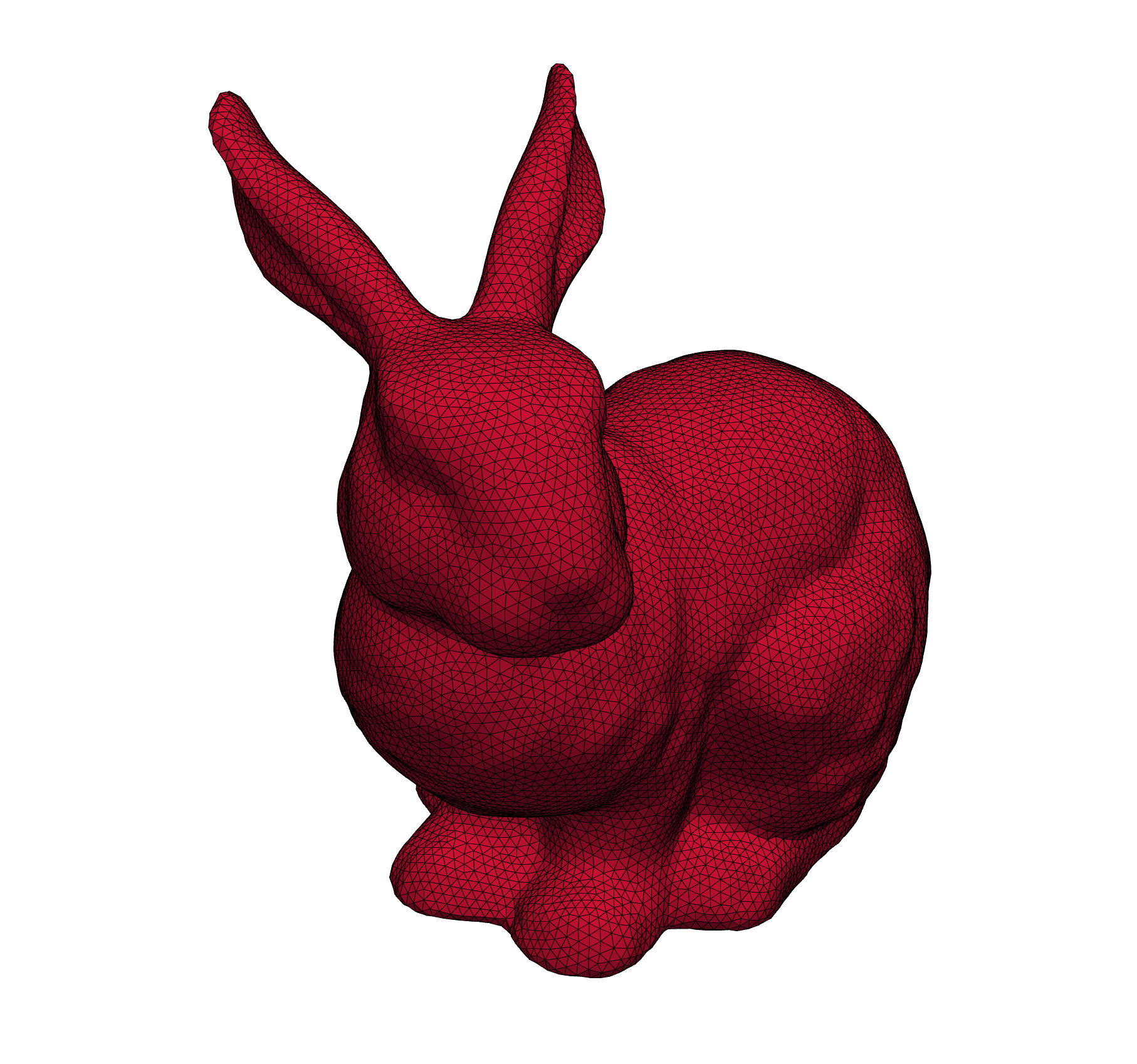}};
\node (bunny_partition) at (3.0,0) {\includegraphics[
width=0.3\textwidth, clip = true, trim = 80 0 80 0]{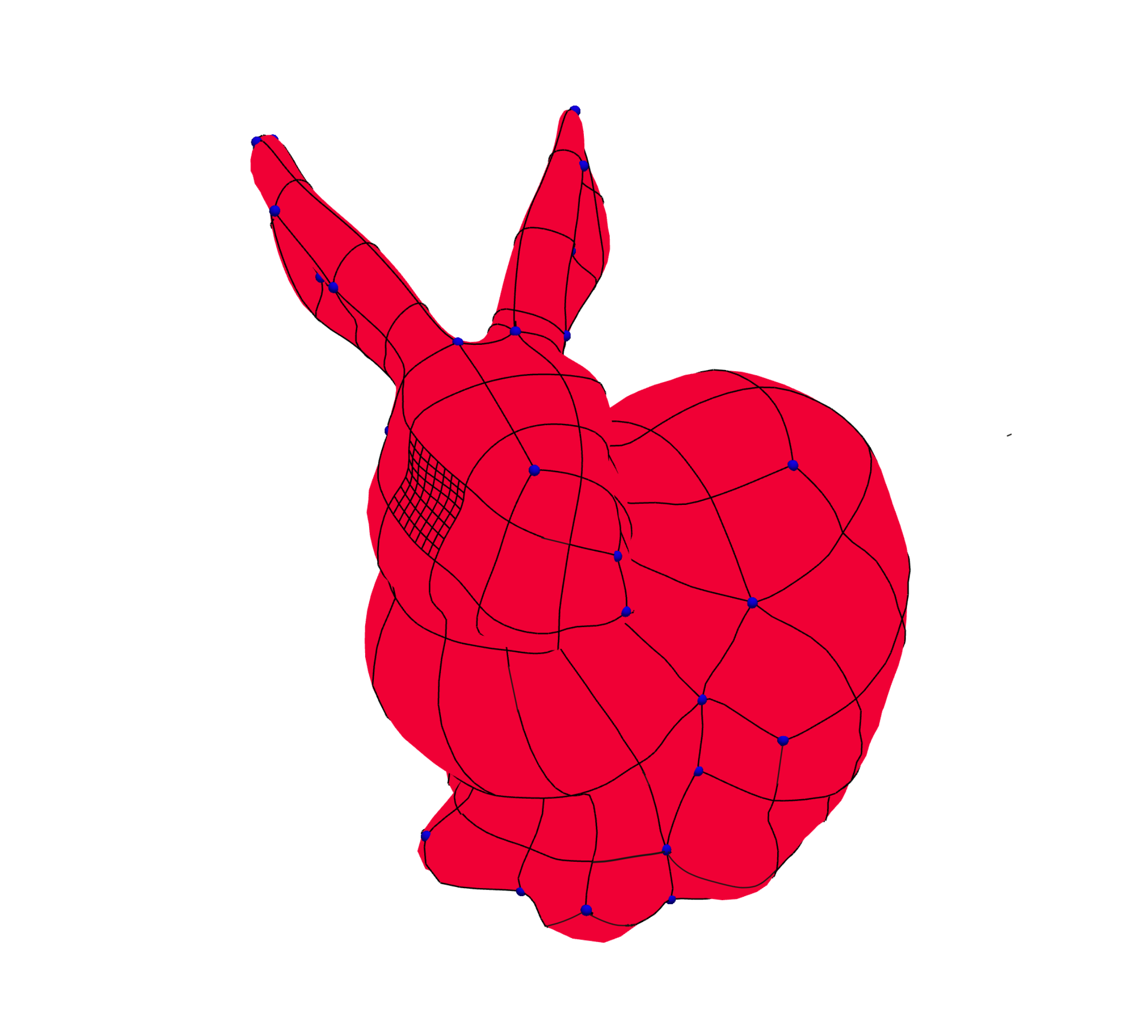}};
\node (a) [] at (-4,2) {};
\node (b) [] at (2.5,2) {};
\draw [ultra thick,black,->] (a) to[bend left] node[above]
(TextNode) {partition} (b);
\end{tikzpicture}
\caption{\label{fig:partition}Partitioning of the Stanford bunny
into 179 quadrilateral patches.}
\end{center}
\end{figure}

The initial quadrangulation for the Stanford bunny
considered in this paper has been taken from \cite{pietroni2016tracing}.
The procedure there is described as follows:
The input is a triangular surface mesh of the geometry endowed with a
\emph{cross field}, which consists of two orthogonal vectors lying
on tangent space of each triangle. The singularities are determined
by localizing points where no local parametrization can be constructed,
see the nodes in Figure~\ref{fig:partition}.
The separatrices are propagated by tracing geodesic paths directly on
the triangular mesh starting from a singularity towards an incident
singularity. By constraining the valence of all singularities to be neither
3 or 5 and avoiding tangential crossings, the separatrices automatically
divide the surface into a conforming quad layout,
see Figure \ref{fig:partition}. Having quadrangulation and patches
described in terms of separatrices at our disposal, we finally fit
a NURBS surface by using the partitioned quadrangulation as control mesh for each patch.
Because the control points match on the edges between any
two incident quad patches, the parametric surface is \(C^0\) globally
and \(C^{p-1}\) patchwise, where \(p\) is the order
of the NURBS representation.

\section{Variational formulation}\label{sec:Varf}
For the discretization of the model problem,
we employ the method of lines. Hence, for fixed \(t>0\),
the discrete variational formulation of the diffusion problem on
random surfaces reads:
\begin{align*}
&\text{find \(u_h({\bs y}) \in V_h({\bs y}) \) such that}\\
&\quad\int_{S({\bs y})}\frac{\partial u_h({\bs y})}
{\partial t}v_h({\bs y})\d{\sigma} 
+\int_{S({\bs y})} \langle\nabla_S u_h({\bs y}),
\nabla_Sv_h({\bs y})\rangle\d\sigma
= \int_{S({\bs y})} f({\bs y})v_h({\bs y}) \d{\sigma}
\end{align*}
for all test functions \(v_h({\bs y}) \in V_h({\bs y})\).
Note that we have omitted here
the dependency of \(u_h\) on \({\bs x}\) and \(t\).
Now, inserting the basis representations of \(u_h({\bs y})\)
and \(v_h({\bs y})\), respectively,
this translates to the linear system 
\[
{\bs M}({\bs y})\frac{\partial}{\partial t} {\bs u}({\bs y},{\bs x},t)
+ {\bs A}({\bs y}){\bs u}({\bs y},{\bs x},t) = {\bs f}({\bs x},t),
\quad {\bs u}({\bs y},{\bs x},0) = {\bs u}_0.
\]
Herein, 
\[
{\bs M}({\bs y}) = [(\varphi_i,\varphi_j)_{L^2(S({\bs y}))}]_{i,j}
\]
denotes the finite element mass matrix, while
\[
{\bs A}({\bs y}) = [(\nabla_{S} \varphi_i,
\nabla_{S} \varphi_j)_{L^2(S({\bs y}))}]_{i,j}
\]
is the finite element stiffness matrix.
Note that the mass matrix and the stiffness matrix are independent
of the time, however, they depend on the random parameter \({\bs y}\).
Since we assume that the initial condition and the source term
are given in material coordinates, their coefficients
\({\bs f}(t)\) and \({\bs u}_0\) are independent of the random
parameter \({\bs y}\). Moreover, we remark that,
due to the definition of the spline basis on the randomly deformed
domain, the coefficient vector \({\bs u}\) in material
and spatial coordinates coincides.

Next, we introduce an appropriate discretization for the time: we employ
the \(\theta\)-scheme
with a uniform time discretization of step size \(\Delta t\).
Thus, we obtain the linear system
\[
({\bs M} + \Delta t\theta {\bs A}){\bs u}(t_{i+1}) = ({\bs M} - \Delta t(1-\theta){\bs A}){\bs u}(t_i) + \Delta t(1-\theta){\bs f}(t_i) + \Delta t \theta {\bs f}(t_{j+1}),
\]
where we have dropped the parameter dependency for a more compact notation.

We close this paragraph by citing the corresponding error estimate.
To this end, we recall
that the surface \(S_0\) consists of \(M\) smooth patches.
A uniform mesh is thus obtained by uniform refinement of the
parameter domain \([0,1]^2\) for each patch and lifting the subdivided
parameter domain to the surface. Hence, subdividing each patch \(j\)
times yields to a mesh of mesh size \(h\sim2^{-j}\).  
For fixed \({\bs y}\), we obtain the following convergence result from \cite{zhu2017isogeometric} under the IGA-\(\theta\) stability
condition introduced therein.

\begin{theorem}\label{thm:conv}
For any fixed \({\bs y}\in[-1,1]^m\), let \(T>0\), \(f\in C^0\left([0,T];L^2\big(S({\bs y})\big)\right), \frac{\partial f}{\partial t}\in L^2\big(S({\bs y})\big),
u(0,{\bs x}) \in H^{p+1}\big(
S({\bs y})\big)\) and the solution \(u\) is such that
\[u\in C^0\left([0,T];H^{p+1}\big(S({\bs y})\big)\right)\text{ and }
\frac{\partial u}{\partial t} \in
L^1\left([0,T];H^{p+1}\big(S({\bs y})\big)\right).
\]
The approximation \(u_h(t_i)\) satisfies \(\frac{\partial u_h}{\partial t}(0) \in L^2\big(S({\bs y})\big)\).
Then, the pointwise \(L^2\)-error for the \(\theta\)-scheme satisfies
\begin{equation}
\begin{aligned}
&\|u(t_i)-u_h(t_i)\|_{L^2(S({\bs y}))}\\
&\qquad\leq C_\theta\bigg[h^{p+1}\bigg(|u(0)|_{H^{p+1}(S({\bs y}))}
+ \left\|\frac{\partial u}{\partial t}\right\|_{L^1([0,T];H^{p+1}(S({\bs y})))}\bigg)\\
&\phantom{leq}\qquad\qquad+\Delta t\bigg(\left\|\frac{\partial u_h}{\partial t}(0)\right\|_{L^2(S({\bs y}))}
+ \left\|\frac{\partial f}{\partial t}\right\|_{L^2([0,T];L^2(S({\bs y})))}\bigg)
\bigg],\quad C_\theta>0.
\end{aligned}
\end{equation}
In the particular case of the Crank-Nicolson method, with additional assumptions \(\frac{\partial^2 f}{\partial t^2}\in L^2\big([0,T]\times S({\bs y})\big)\) and \(\frac{\partial^2 u_h}{\partial t^2} \in L^2\big(S({\bs y})\big)\) we obtain
\begin{equation}
\begin{aligned}
&\|u(t_i)-u_h(t_i)\|_{L^2(S({\bs y}))}\\
&\qquad\leq{C}_{1/2}\bigg[h^{p+1}\bigg(|u(0)|_{H^{p+1}(S({\bs y}))}
+\left\|\frac{\partial u}{\partial t}\right\|_{L^1([0,T];H^{p+1}(S({\bs y})))}\bigg)\\
&\qquad\qquad\quad\phantom{\leq}+(\Delta t)^2\bigg(\left\|\frac{\partial^2 u_h}{\partial t^2}(0)\right\|_{L^2(S({\bs y}))}
+ \left\|\frac{\partial^2 f}{\partial t^2}\right\|_{L^2([0,T];L^2(S({\bs y})))}\bigg)
\bigg],
\end{aligned}
\end{equation}
where \(u(t_i)\) refers to the exact solution at time \(t_i\),
\(u_h\in V_h({\bs y})\) is the Galerkin approximation
and \(p\) is the minimum of \(p_1\) and \(p_2\),
cp.\ \eqref{eq:AnsatzSpace}. Herein, \(C_\theta\) and \(C_{1/2}\) are independent of \(\Delta t\) and \(h\).
\end{theorem}

\section{Computation of quantities of interest}\label{sec:qoi}
\subsection{Quantities of interest}
We consider quantities of interest
that are defined with respect to the unperturbed
reference surface \(\Sref\).
There holds for an arbitrary time dependent function
\(u_h(t) \in V_h({\bs y})\) that
\begin{align*}
\hat{u}_h({\bs y},\xref,t)&=
\big(u_h(t) \circ\DefField\big)({\bs y},\xref)=\sum_{i=1}^{N_{{\bs p},{\bs\Xi}}}
u_i({\bs y},t)(\varphi_i\circ\DefField)({\bs y},\xref)\\
&=\sum_{i=1}^{N_{{\bs p},{\bs\Xi}}}
u_i({\bs y},t)\hat{\varphi}_i(\xref)=\widehat{\bs\Phi}{\bs u}({\bs y},t)
\end{align*}
with the time and parameter dependent coefficient vector
\[{\bs u}({\bs y},t)
\isdef[u_1({\bs y},t),u_2({\bs y},t),\ldots,u_{N_{{\bs p},{\bs\Xi}}}({\bs y},t)]^\intercal.
\]
This means that, by construction, the coefficients of the basis
representations in \(\Scal_{{\bs p},{\bs\Xi}}\big(S({\bs y})\big)\)
and \(\Scal_{{\bs p},{\bs\Xi}}(\Sref)\) coincide.
Consequently, we obtain the following two equations for the expectation and correlation respectively,
\[
\Ebb[\hat{u}_h(t)](\xref)
=\int_\Omega \widehat{\bs\Phi}(\xref){\bs u}({\bs y},t)\d\Pbb(\omega)
=\widehat{\bs\Phi}(\xref)\int_\Omega {\bs u}({\bs y},t)\d\Pbb(\omega)
=\widehat{\bs\Phi}(\xref)\Ebb[{\bs u}(t)],
\]
and 
\[
\Cor[\hat{u}_h]\big((\xref,t),(\xref',t')\big)=\widehat{\bs\Phi}(\xref)
\Ebb[{\bs u}(t){\bs u}^\intercal(t')]
\widehat{\bs\Phi}^\intercal(\xref').
\]
Hence, the quantities of interest under consideration can be
approximated by computing the corresponding quantities of interest
of the coefficient vectors, i.e.,
\(\Ebb[{\bs u}(t)]\) and \(\Ebb[{\bs u}(t){\bs u}^\intercal(t')]\).
Numerically, we approximate those quantities by applying suitable
quadrature formulas in the parameter, such as the Monte Carlo
method and the quasi-Monte Carlo quadrature
based on Halton points, cf.\ \cite{Caf98}. Thus,
we end up with approximations
\begin{equation}\label{eq:expectation}
\Ebb[\hat{u}_h(t)]
\approx \widehat{\bs\Phi}
\bigg(\frac{1}{N_q}\sum_{i=1}^{N_q} {\bs u}({\bs\xi}_i,t)\bigg),
\end{equation}
and 
\begin{equation}
\Cor[\hat{u}_h]\big((\xref,t),(\xref',t')\big)
\approx \widehat{\bs\Phi}(\xref)\bigg(\frac{1}{N_q}
\sum_{i=1}^{N_q}{\bs u}({\bs\xi}_i,t) {\bs u}({\bs\xi}_i,t')^\intercal
\bigg) \widehat{\bs\Phi}^\intercal(\xref'),
\end{equation}
where \({\bs\xi}_i\in\Gamma_m\isdef[-1,1]^m\), \(i=1,\ldots, N_q\),
are the sample points.
Quantities of interest based on linear functionals
of the solution can be dealt with
in a similar fashion. In summary, we remark that the
computation of quantities of interest amounts to
high dimensional quadrature problems for the coefficient
vector of the basis representation in 
\(\Scal_{{\bs p},{\bs\Xi}}(\Sref)\). 

For the computation of the expectation,
we only need to store one coefficient vector for each desired time point.
In total, the required memory size is thus \(\mathcal{O}(N_tN_{{\bs p},{\bs\Xi}})\).
In the case of the correlation, however, the situation is much worse,
as we need to store an \((N_tN_{{\bs p},{\bs\Xi}}) \times (N_tN_{{\bs p},{\bs\Xi}})\)
matrix, resulting in a memory size of \(\mathcal{O}(N_t^2N_{{\bs p},{\bs\Xi}}^2)\).
This matrix easily exceeds the memory size of a large computer
if \(N_t\) or \(N_{{\bs p},{\bs\Xi}}\) is large. This motivates
the low rank representation of the correlation, which is discussed
in the subsequent paragraph.

\subsection{Low rank approximation of the correlation}
For the cost efficient storage of the correlation, we employ an online
algorithm that keeps track of the most important directions for the storage
of the correlation matrix. To this end, we consider the online singular value
decomposition (SVD) suggested in \cite{brand2002incremental}.
Given the sequence of time-dependent coefficient vectors
\(\{{\bs u}({\bs\xi}_i,t)\}_{i=1}^{N_q}\) defined by the evaluation at
the sample points, we employ the online singular value decomposition
for each discrete time step. To outline the algorithm, we define the
matrices
\[
{\bs W}_\ell(t_i)\isdef
[{\bs u}^{(1)}(t_i),\ldots,
{\bs u}^{(\ell)}(t_i)]\isdef
[{\bs u}({\bs\xi}_1,t_i),\ldots,
{\bs u}({\bs\xi}_{\ell},t_i)]\in\Rbb^{N_{{\bs p},{\bs\Xi}}\times\ell},
\]
which contain the first \(\ell\) samples at time \(t=t_i\).
The coefficient matrix of the correlation of the solution between
the times \(t_i\) and \(t_{i'}\) for the first \(\ell\)
samples is thus given by
\[
{\bs C}_\ell(t_i,t_{i'})=\frac{1}{\ell}{\bs W}_\ell(t_i)
{\bs W}_\ell^\intercal(t_{i'}).
\]
We denote the SVD of \({\bs W}_\ell(t_i)\) by
\[
{\bs W}_\ell(t_i) = {\bs U}_\ell(t_i){\bs \Sigma}_\ell(t_i)
{\bs V}_\ell^\intercal(t_i).
\]
In the \(\ell\)-th step, the vector \({\bs u}^{(\ell)}(t_i)\)
is used to produce a potential new column vector of the matrix
\({\bs U}_\ell(t_i)\). 
Accounting for the situation of limited memory, we only consider
the first \(k\leq N_q\) most important columns of \({\bs U}_\ell(t_i)\).
Here, the importance is determined by the magnitude of the associated
singular value stored in the matrix \({\bs \Sigma}_\ell(t_i)\). As we 
only wish to compute the correlation, we do not need to explicitly
store the right-orthonormal matrices \({\bs V}_\ell(t_i)\).
Instead, we rather store a low dimensional coefficient matrix
\(\lowC_\ell(t_i,t_{i'})
\isdef \frac{1}{\ell}{\bs V}_\ell^\intercal(t_i) {\bs V}_\ell(t_{i'})\)
between the times \(t_i\) and \(t_{i'}\). Given that we only keep track of the first
\(k\) most important vectors, the size of each \(\lowC_\ell(t_i,t_{i'})\)
is never larger than \(k \times k\).

In what follows, we will discuss how to compute the
low rank approximation of the space-time coefficient matrix in detail
by adapting the idea from \cite{brand2002incremental}.
Given a new sample \({\bs u}^{(\ell)}(t_i)\) for \(i=1,\cdots,N_t\),
the first step is computing the coefficient vector and the residual
at each discrete time \(t_i\) by projecting the new sample solution
onto the subspace spanned by the orthonormal basis
\({\bs U}_{\ell-1}(t_i)\). This way, we obtain the coefficients
\[
{\bs c}(t_i) = {\bs U}_{\ell-1}^\intercal(t_i) {\bs u}^{(\ell)}(t_i)
\]
and the remainder
\[
\hat{\bs r}(t_i) = {\bs u}^{(\ell)}(t_i) - {\bs U}_{\ell-1}(t_i) {\bs c}(t_i).
\]
The associated normalized remainder is given by
\[
{\bs r}(t_i)=\frac{\hat{\bs r}(t_i)}{r(t_i)}\quad
\text{with }r(t_i)\isdef\|\hat{\bs r}(t_i)\|_2.
\]
Next, we pad \({\bs U}_{\ell-1}(t_i)\) by \({\bf r}(t_i)\)
and \({\bs \Sigma}_{\ell-1}(t_i)\)
by adding \({\bs c}(t_i)\) and \(r(t_i)\) according to
\[
\hat{\bs U}_\ell(t_i) \isdef \begin{bmatrix}
{\bs U}_{\ell-1}(t_i) & {\bs r}(t_i)\\
\end{bmatrix}
\]
and 
\[
\hat{\bs \Sigma}_\ell(t_i) \isdef \begin{bmatrix}
{\bs \Sigma}_{\ell-1}(t_i) & {\bs c}(t_i)\\
{\bs 0} & r(t_i)
\end{bmatrix}.
\]
Especially, we have \(
{\bs W}_\ell(t_i)= \hat{\bs U}_\ell(t_i)\hat{\bs \Sigma}_\ell(t_i){\hat{\bs V}_\ell}^\intercal(t_i),\) where
\[
\hat{\bs V}_\ell(t_i) \isdef \begin{bmatrix}
{\bs V}_{\ell-1}(t_i) & {\bs 0}\\
{\bs 0} & 1
\end{bmatrix}.
\]
In order to store the low rank approximation of the
space-time coefficient matrix, we need to update
the small matrix \(\lowC_\ell(t_i,t_{i'})\) between \(t_i\) and \(t_{i'}\)
 for \(i,i' = 1,\cdots,N_t\). To this end, we pad it according to
\[
\hat{\lowC}_\ell(t_i,t_{i'}) \isdef
\frac{1}{\ell}{\hat{\bs V}_\ell}^\intercal(t_i)
\hat{\bs V}_\ell(t_{i'})=\frac{\ell-1}{\ell}\begin{bmatrix}
\lowC_{\ell-1}(t_i,t_{i'}) & {\bs 0}\\
{\bs 0} & \frac{1}{\ell-1}
\end{bmatrix}.
\]
In the final step, we diagonalize \(\hat{\bs \Sigma}_\ell(t_i)\),
which amounts to the SVD of a low dimensional matrix according to
\(\hat{\bs \Sigma}_\ell(t_i) = {\bs P}_\ell(t_i) {\bs \Sigma}_\ell(t_n) {{\bs Q}_\ell}^\intercal(t_i)\) and perform the updates
\begin{equation}\label{eq:updateU}
{\bs U}_\ell(t_i) =\hat{\bs U}_\ell(t_i){\bs P}_\ell(t_i),
\end{equation}
and
\[
\lowC_\ell(t_i,t_{i'}) 
= {{\bs Q}_\ell}^\intercal(t_i)\hat{\lowC}_\ell(t_i,t_{i'}) {{\bs Q}_\ell}(t_{i'}).
\]

If the rank of \({\bs U}_\ell(t_i)\) is larger than the maximally allowed number \(k\),
we truncate \({\bs U}_\ell(t_i)\) by choosing the left most \(k\) columns, given that
the singular values are stored in decreasing order. Moreover, we truncate
\({\bs \Sigma}_\ell(t_i)\) by keeping the top-left \(k\times k\) block.
Finally, we truncate \(\lowC_\ell(t_i,t_{i'})\) by keeping the top most \(k\) rows if
the rank of \({\bs U}_\ell(t_i)\) is larger than \(k\) and truncate
\(\lowC_\ell(t_i,t_{i'})\) by keeping the left most \(k\) columns if the rank of \({\bs U}_\ell(t_{i'})\) is larger than \(k\). The procedure is summarized in
Algorithm \ref{alg:lowrank}.

\begin{algorithm}[htb]
	\KwIn{stream of solutions \(\{{\bs u}(t,{\bs\xi}_i)\}_{i=1}^{N_q}\)} 
	\KwOut{${\bs U}(t_i)$, ${\bs \Sigma}(t_i)$ and $\lowC(t_i,t_{i'})$} 
	\For{$j=1$ \KwTo $N_q$}{ 
			\For{$i= 1$ \KwTo $N_t$}{
	    $\hat{\bs r}(t_i) = {\bs u}(t_i,{\bs\xi}_i)$;\\
		\uIf{$j=1$}{
			$[{\bs U}(t_i), {\bs \Sigma} (t_i), {\bs V}(t_i)] = \operatorname{SVD}(\hat{\bs r}(t_i))$;\\
		}
		\Else{
		    ${\bs c}(t_i) = {\bs U}^\intercal(t_i)\hat{\bs r}(t_i)$\\
		    $\hat{\bs r}(t_i) = {\bs u}_i(t_i) - {\bs U}(t_i){\bs c}(t_i)$\\
		    $r(t_i)=\|\hat{\bs r}(t_i)\|_2$\\
		    ${\bs r}(t_i) = \hat{\bs r}(t_i)/r(t_i)$\\
		    $\hat{\bs U}(t_i) = [{\bs U}(t_i),{\bs r}(t_i)]$;\\
		    $\hat{\bs \Sigma}(t_i) = [{\bs \Sigma}(t_i),{\bs c}(t_i);{\bs 0},r(t_i)]$;\\
		    $[{\bs P}(t_i), {\bs \Sigma} (t_i), {\bs Q} (t_i)] = \operatorname{SVD}(\hat{\bs \Sigma} (t_i))$;\\
		    ${\bs U}(t_i) = \hat{\bs U}(t_i){\bs P} (t_i)$;
		}
		}
		\For{$i =1$ \KwTo $N_t$}{
		\For{$i' =i$ \KwTo $N_t$}{
		    \uIf{$j=1$}{
		    $\lowC(t_i,t_{i'}) = {{\bs V}}^{\intercal} (t_i){\bs V}(t_{i'})$;\\
		    }
		    \Else{
		    $\hat{\lowC}(t_i,t_{i'}) =\frac{j-1}{j}[\lowC(t_i,t_{i'}),{\bs 0};{\bs 0},\frac{1}{j-1}]$;\\
		    $\lowC(t_i,t_{i'}) = {{\bs Q}}^\intercal (t_i) \hat{\lowC}(t_i,t_{i'}) {{\bs Q}} (t_{i'})$;
		    }
		\If{$j > k$}{
		    $\lowC(t_i,t_{i'}) = \lowC(t_i,t_{i'})(1:k,1:k)$
		    }
		}
				    \If{$j > k$}{
		    ${\bs U}(t_i) = {\bs U}(t_i)(:,1:k)$\\
		    ${\bs \Sigma}(t_i) = {\bs \Sigma}(t_i)(1:k,1:k)$\\
		    }
		}
	}
	\caption{Low rank approximation of the space time correlation} 
	\label{alg:lowrank} 
\end{algorithm}

We arrive at a tall-and-skinny matrix \({\bs U}_\ell(t_i)\) 
whose columns are used as low rank basis for time \(t_i\).
In particular, the matrix of singular values \({\bs \Sigma}_\ell(t_i)\) 
gives us a means to track the importance of each
basis vector. Moreover, the matrix \(\lowC_\ell(t_i,t_{i'})\) represents
the correlation structure between time points
 \(t_i\) and \(t_{i'}\) in a compressed manner. In summary, we end up
 with the low rank approximation
\[
{\bs C}_{N_q}(t_i,t_{i'}) \approx \tilde{{\bs C}}_{N_q}(t_i,t_{i'}) \isdef {\bs U}_{N_q}(t_i) {\bs \Sigma}_{N_q}(t_i) \lowC_{N_q}(t_i,t_{i'}) {\bs \Sigma}_{N_q}^\intercal(t_{i'}) {\bs U}_{N_q}^\intercal(t_{i'}).
\]
The corresponding low rank space time correlation is hence given by
\begin{align*}
\Cor[\hat{u}_h]\big((\xref,t_i),(\xref',t_{i'})\big)
\approx\widehat{\bs\Phi}(\xref)
\tilde{{\bs C}}_{N_q}(t_i,t_{i'})
\widehat{\bs\Phi}^\intercal(\xref').
\end{align*}

We conclude this paragraph by discussing the computational cost
of the low rank representation. Storing the matrices \(\lowC(t_i,t_{i'})\),
\(i,i'=1,\ldots,N_t\) results in cost of \(\mathcal{O}(N_tk^2)\), while
the cost for storing the bases \({\bs U}(t_i)\) is of cost
\(\mathcal{O}(N_tN_{{\bs p},{\bs\Xi}}k)\). Assuming 
\(k\leq N_{{\bs p},{\bs\Xi}}\), the storage cost is
hence of order \(\mathcal{O}(N_tN_{{\bs p},{\bs\Xi}}k)\), which is a huge
reduction compared to the original cost of
\(\mathcal{O}(N_t^2N_{{\bs p},{\bs\Xi}}^2)\).
The computational cost for obtaining the low rank approximation is
comprised of the projection steps of cost
\(\mathcal{O}(N_qN_tN_{{\bs p},{\bs\Xi}}k)\),
while the cost of the required SVD in each step is
\(\mathcal{O}(N_qN_t^2k^3)\) and the cost of the update formula
\eqref{eq:updateU} is \(\mathcal{O}(N_tN_{{\bs p},{\bs\Xi}}N_qk^2)\).
This amounts to an overall computation cost of
\(\mathcal{O}(N_qN_tN_{{\bs p},{\bs\Xi}}k^2+N_qN_t^2k^3)\), which for moderate
\(k\) is again much smaller than the original cost of 
\(\mathcal{O}\big(N_q(N_tN_{{\bs p},{\bs\Xi}})^2\big)\).
\section{Numerical experiments}\label{sec:Numerics}
\subsection{Convergence of the spatial solver}
All computations for this article have been performed with the 
Boundary Element Method Based Engineering Library (\verb+Bembel+), 
see \cite{dolz2020bembel}. In particular, we have extended \verb+Bembel+
to also support local operators, like the Laplace-Beltrami operator, and
added a class for the fast deformation of computational geometries.
We numerically validate our implementation for the diffusion equation
by showing that the theoretical convergence rates from Theorem~\ref{thm:conv}
are obtained. At first, we consider a diffusion equation on the unit sphere, 
i.e.,
\begin{align}\label{eq:heatex_sim}
\begin{cases}
\partial_t u({\bs x},t)
- \LBref u({\bs x},t) = 0, &\quad {\bs x}
\in \mathbb{S}^2,\\
u({\bs x},0) = x_3, &\quad {\bs x} \in \mathbb{S}^2.\\
\end{cases}
\end{align}

\begin{figure}[htb]
\begin{center}
\pgfplotsset{width=0.5\textwidth}
\pgfplotsset{minor grid style={dotted,black!30}}
\pgfplotsset{grid style={solid,black!30}}
\begin{tikzpicture}
\begin{semilogyaxis}[ytick distance=10,grid=both, ymin= 1e-9, ymax = 0.1, xmin = 0, xmax = 4, 
legend style={legend pos=south west,font=\small}, legend cell align={left},%
ylabel={\(L^2\)-error}, xlabel ={level \(j\)}]
\addplot[line width=1pt,color=red,mark=oplus] table[x index=0,y index=1]{./Data/sphere_l2_error_1.txt};\addlegendentry{$p=1$};
\addplot[line width=1pt,color=green!70!black,mark=oplus] table[x index=0,y index=1]{./Data/sphere_l2_error_2.txt};\addlegendentry{$p=2$};
			\addplot[line width=1pt,color=blue,mark=oplus] table[x index=0,y index=1]{./Data/sphere_l2_error_3.txt};\addlegendentry{$p=3$};
						\addplot[dashed] table[x index=0,y expr=0.08*2^-2*\thisrowno{0}]{./Data/sphere_l2_error_1.txt};
\addplot[dashed] table[x index=0,y expr=0.08*2^-3*\thisrowno{0}]{./Data/sphere_l2_error_1.txt};
\addplot[dashed] table[x index=0,y expr=0.08*2^-4*\thisrowno{0}]{./Data/sphere_l2_error_1.txt};
						\addlegendentry{$h^2,h^3,h^4$};
			\end{semilogyaxis}
			\end{tikzpicture}
	\end{center}
	\caption{\(L^2\)-error versus the level \(j\) on the unit sphere.}
	\label{fig:sphere}
\end{figure}
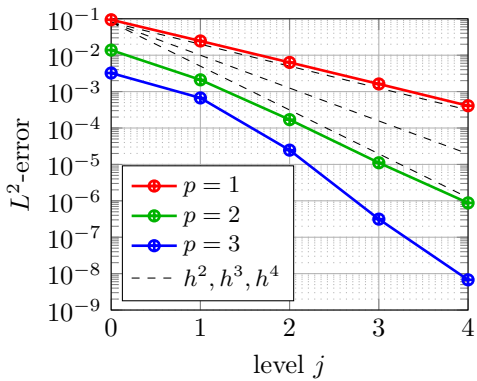

The corresponding exact solution is
\(u({\bs x},t)=e^{-2t}x_3\).
Figure~\ref{fig:sphere} shows the corresponding convergence plot for 
a representation of the sphere with \(M=6\) patches and a uniform
refinement up to level \(4\) and \(\Delta t=10^{-5}\),
evaluated for the final time \(T=1\). For the time stepping we 
use the Crank-Nicolson method and we consider \(p=1,2,3\).
It can be seen that the theoretical rate from Theorem~\ref{thm:conv}
is achieved. Indeed,
 the observed rate appears to be even higher.
 
\begin{figure}[htb]
\begin{center}
			\includegraphics[scale=0.08,clip=true,trim=2000 350
			2000 400]{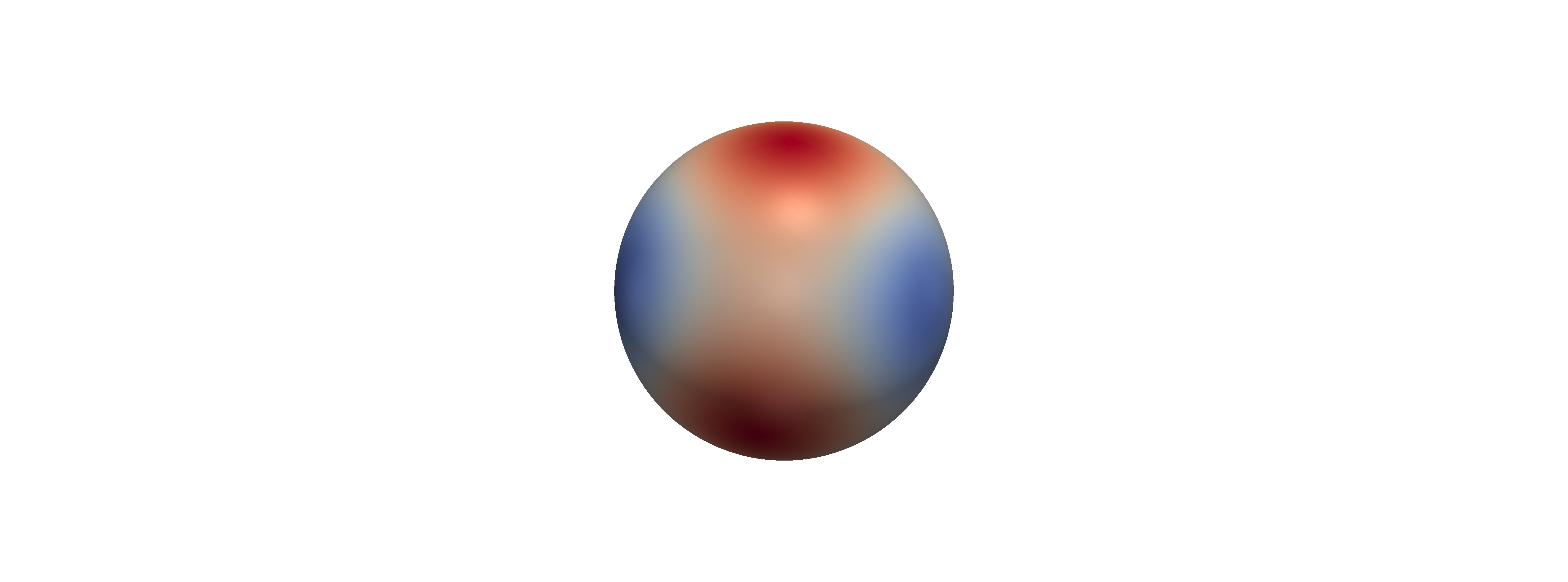}
			\includegraphics[scale=0.06,clip=true,trim=1800 200
			1800 100]{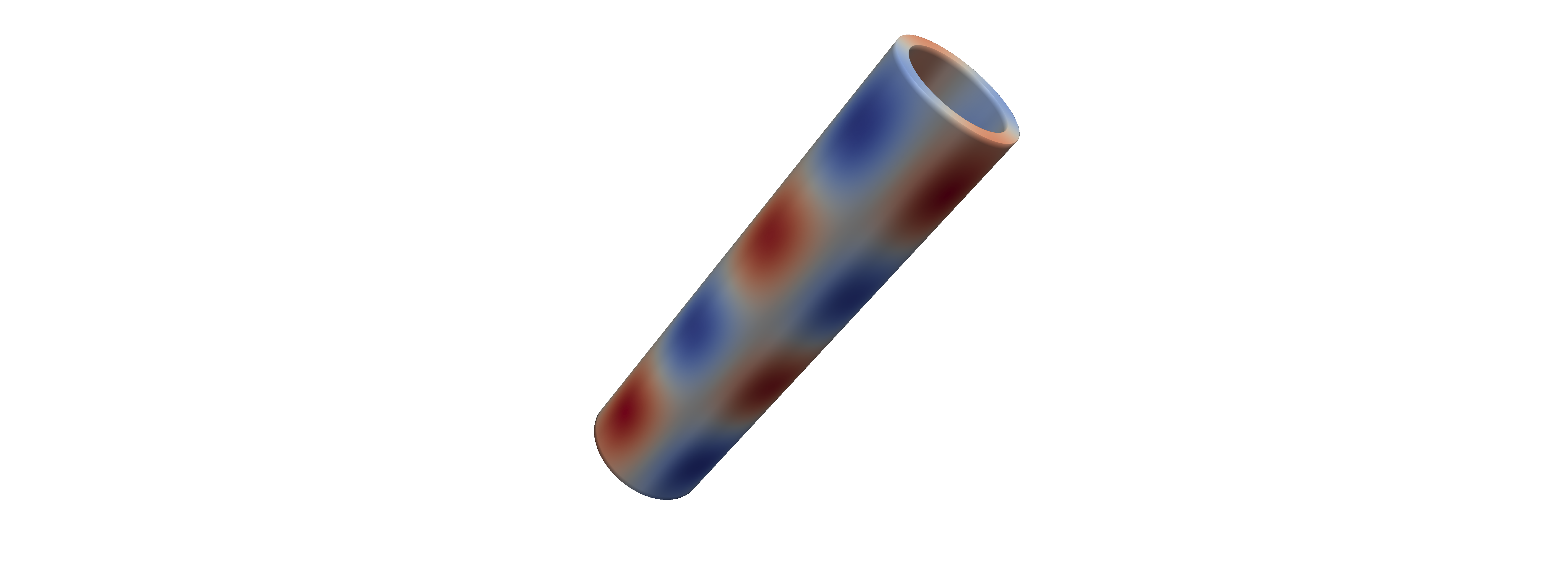}\qquad
			\includegraphics[scale=0.07,clip=true,trim=2000 200
			2000 350]{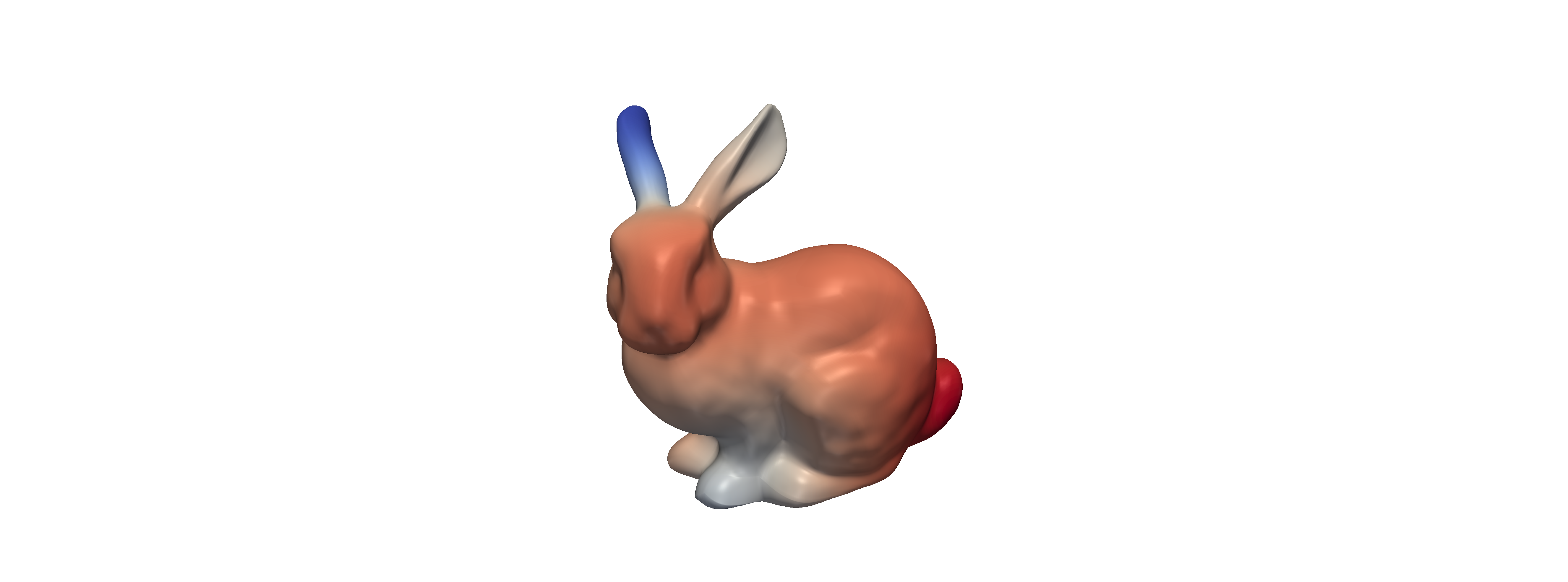}
	\caption{\label{fig:diffgeos}Visualization of the solutions at time \(T=1\).}
	\end{center}
\end{figure}

Next, in view of the previous numerical results
and Theorem~\ref{thm:conv}, we fix \(p=2\) and
\(\Delta t = 10^{-3}\), such that the 
convergence is not limited by the time discretization. We consider the
diffusion problem 
\begin{align*}
\begin{cases}
\partial_t u({\bs x},t)
- \LBref u({\bs x},t) = \sin(\pi x_1)\sin(\pi x_2)\sin(\pi x_3), &\quad {\bs x}
\in S,\\
u({\bs x},0) = x_3, &\quad {\bs x} \in S.\\
\end{cases}\end{align*}

for three different geometries. Namely, these are the unit sphere
(bounding box: \([-1,1]^3\)),
a pipe geometry (bounding box: \([0,2.4]\times[-0.3,0.3]^2)\) and the Stanford bunny
(bounding box: \([-0.95,0.61]\times[0.33,1.86]\times[-0.61,0.59])\).
Vizualisations of the solutions for
\(T=1\) can be found in Figure~\ref{fig:diffgeos}.

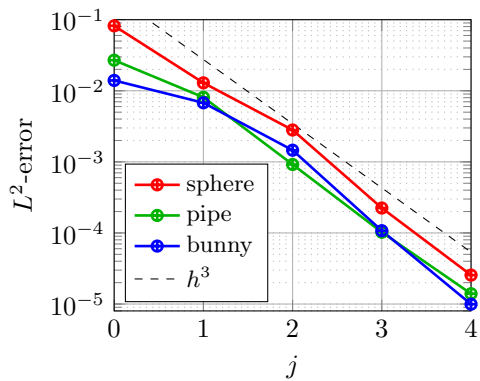
\begin{figure}[htb]
\begin{center}
\pgfplotsset{width=0.5\textwidth}
\pgfplotsset{minor grid style={dotted,black!30}}
\pgfplotsset{grid style={solid,black!30}}
\begin{tikzpicture}
\begin{semilogyaxis}[grid=both, ymin= 0.8e-5, ymax = 0.1, xmin = 0, xmax = 4,%
	   legend style={legend pos=south west,font=\small}, legend cell align={left},%
	    ylabel={\(L^2\)-error}, xlabel ={\(j\)}]
\addplot[line width=1pt,color=red,mark=oplus] table[x index=0,y index=1]{./Data/l2_error.txt};\addlegendentry{sphere};
\addplot[line width=1pt,color=green!70!black,mark=oplus] table[x index=0,y index=2]{./Data/l2_error.txt};\addlegendentry{pipe};
\addplot[line width=1pt,color=blue,mark=oplus] table[x index=0,y index=3]{./Data/l2_error.txt};\addlegendentry{bunny};
\addplot[dashed] table[x index=0,y expr=0.22*2^-3*\thisrowno{0}]{./Data/l2_error.txt};\addlegendentry{$h^3$};
\end{semilogyaxis}
\end{tikzpicture}
\caption{\label{fig:geoconv}Average \(L^2\)-error decrease for \(p=2\) and the different geometries under consideration.}
\end{center}
\end{figure}

For the convergence study, we compute the average \(L^2\)-error evaluated 
at the time grid \(t=0.001,0.002,\ldots,1\), i.e.,
\[
e_{L^2}^{(j)} = \frac{1}{N_t}\sum_{i=1}^{N_t} \norm{u_J(t_i)-u_j(t_i)}_{L^2(S)},\quad
N_t=1000,
\]
where we use the numerical solution on level \(J=5\) as ground truth.
The corresponding error plot can be found in
Figure~\ref{fig:geoconv}. 
As can be seen, the error is almost divided by \(8\) when the level \(j\) increases by one,
which reflects the expected convergence rate of \(h_j\sim 2^{-3j}\).

\subsection{Convergence of the low rank approximation}
In this part, we numerically test the convergence of the proposed algorithm 
with respect to the predefined rank \(k\), the generalization error
and the stability of the algorithm with respect to the magnitude of the domain
perturbation. Now we consider the diffusion problem
\begin{align*}
\begin{cases}
\partial_t u({\bs y},{\bs x},t)
- \LBref u({\bs y},{\bs x},t) = \sin(\pi x_1)\sin(\pi x_2)\sin(\pi x_3), &\quad {\bs x}
\in S({\bs y}),\\
u({\bs y},{\bs x},0) = x_3, &\quad {\bs x} \in S({\bs y}).\\
\end{cases}
\end{align*}

As before, we consider the randomly deformed sphere, pipe and 
Stanford bunny.
For the time discretization, we use again the Crank-Nicolson method
with time steps \(\Delta t = 0.001\) for the sphere and the pipe and  
\(\Delta t = 0.01\) for the Stanford bunny.
For the spatial discretization, we set the level and the polynomial 
order of the basis functions with \(j=4\) for the sphere and the pipe
and \(j=3\) for the Stanford bunny. Moreover, we always choose the
polynomial degree \(p=2\).
For the sphere and pipe, this setting results in an
accuracy of \(10^{-5}\), see Figure~\ref{fig:geoconv}. 
For the Stanford bunny, the computational cost is much larger
than for the sphere and pipe because of the more complex geometry.
Therefore, we reduce the level \(j\) by one, resulting in an
accuracy of \(10^{-4}\), see Figure~\ref{fig:geoconv}.
In all cases, the choice of the time step ensures that the
error is not dominated by the time discretization.
This leads to \(N_{{\bs p},{\bs\Xi}} = 1736, 6936, \text{ and } 14501\)
for the sphere, pipe and the Stanford bunny, respectively,
and results in a spatial discretization
error about \(10^{-4}\) for all (undeformed) geometries, cp.\ previous
example.

\begin{figure}[htb]
\begin{center}
\includegraphics[width=0.3\textwidth,trim=700pt 100pt 700pt 100pt,clip=true]{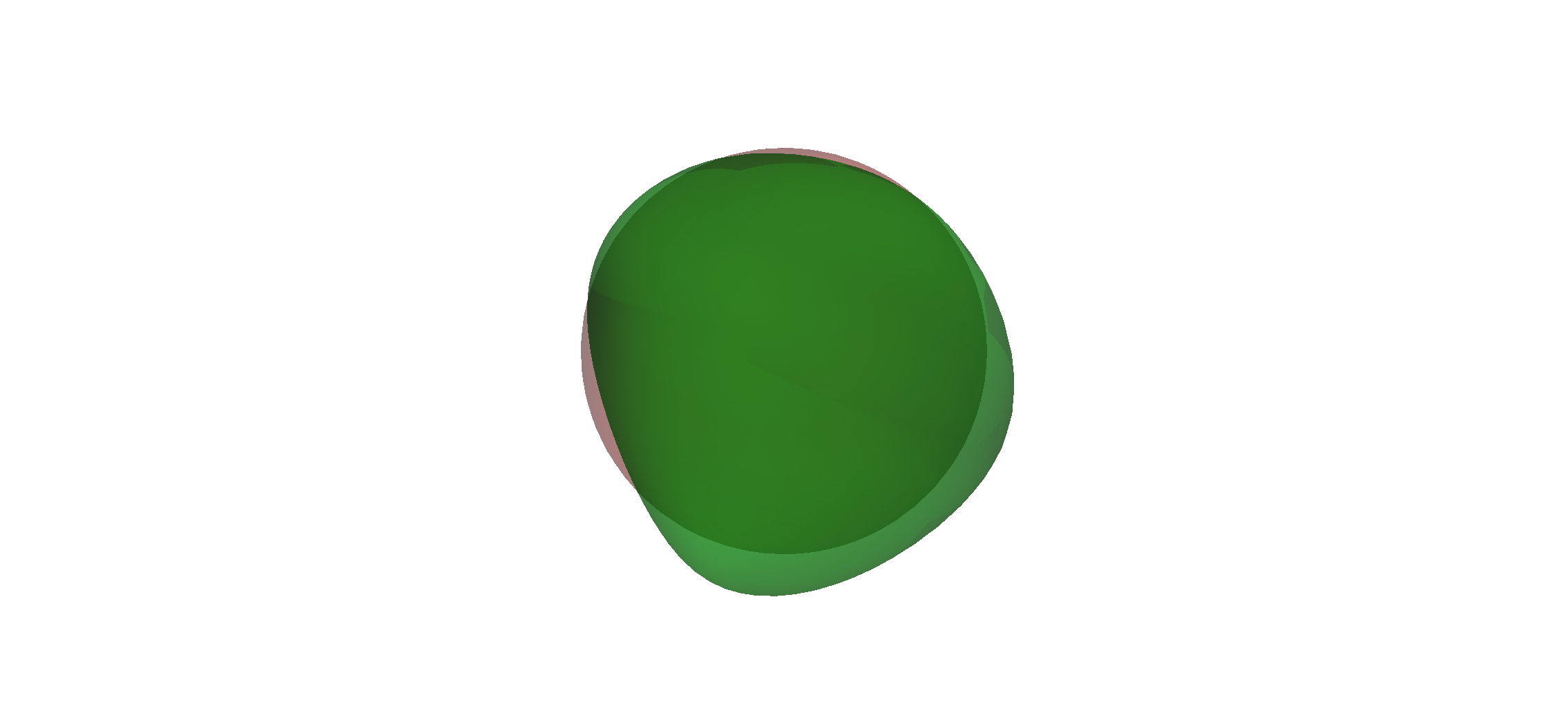}
\includegraphics[width=0.3\textwidth,trim=700pt 100pt 700pt 100pt,clip=true]{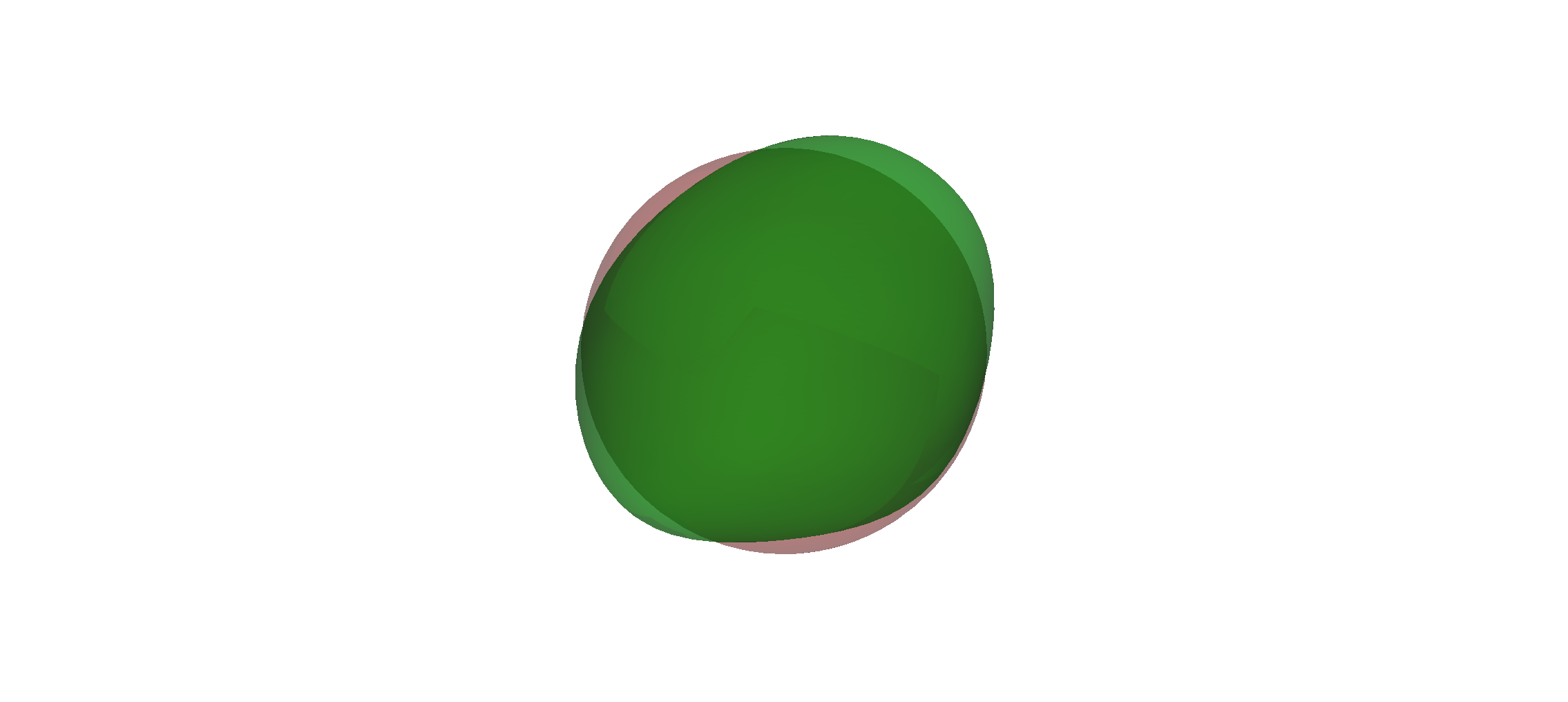}
\includegraphics[width=0.3\textwidth,trim=700pt 100pt 700pt 100pt,clip=true]{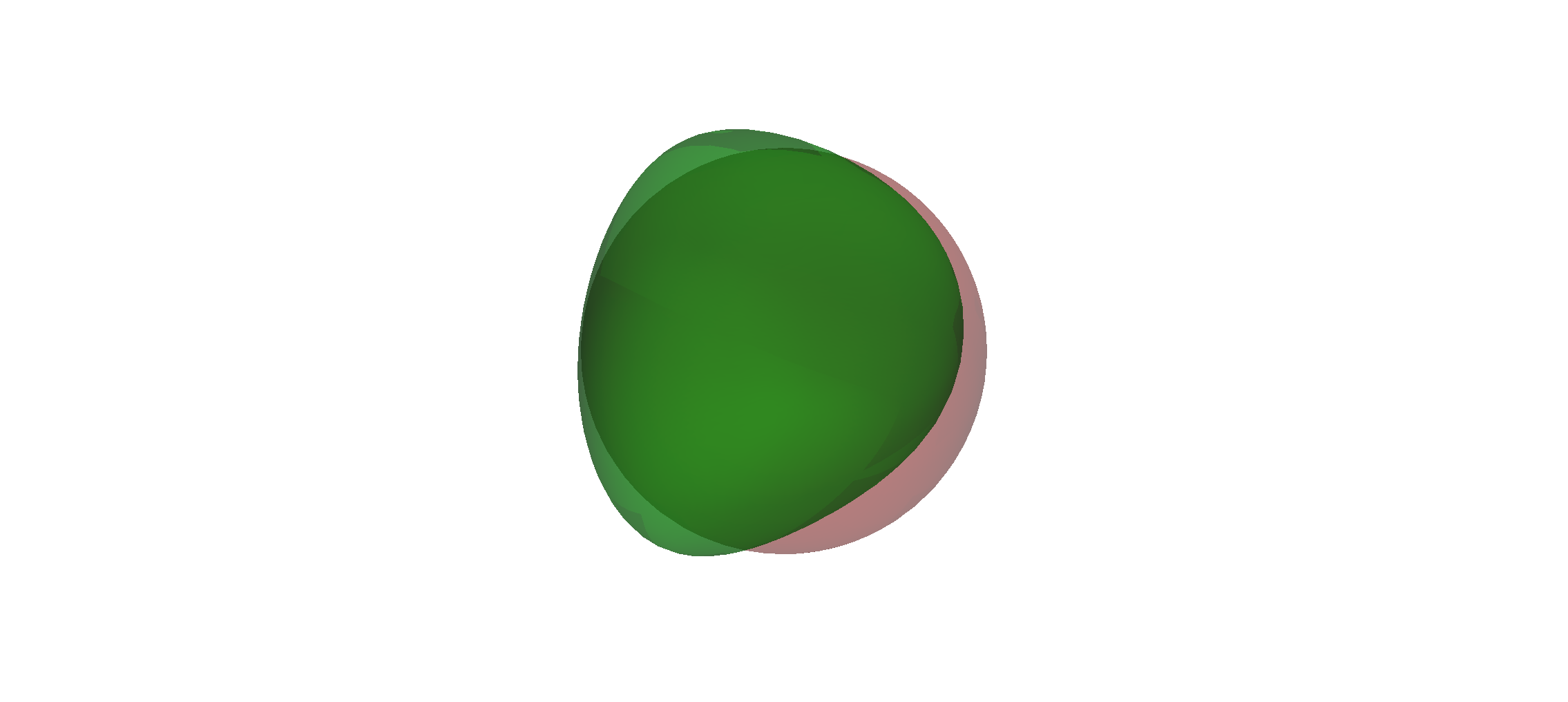}
\includegraphics[width=0.3\textwidth,trim=700pt 50pt 550pt 0pt,clip=true]{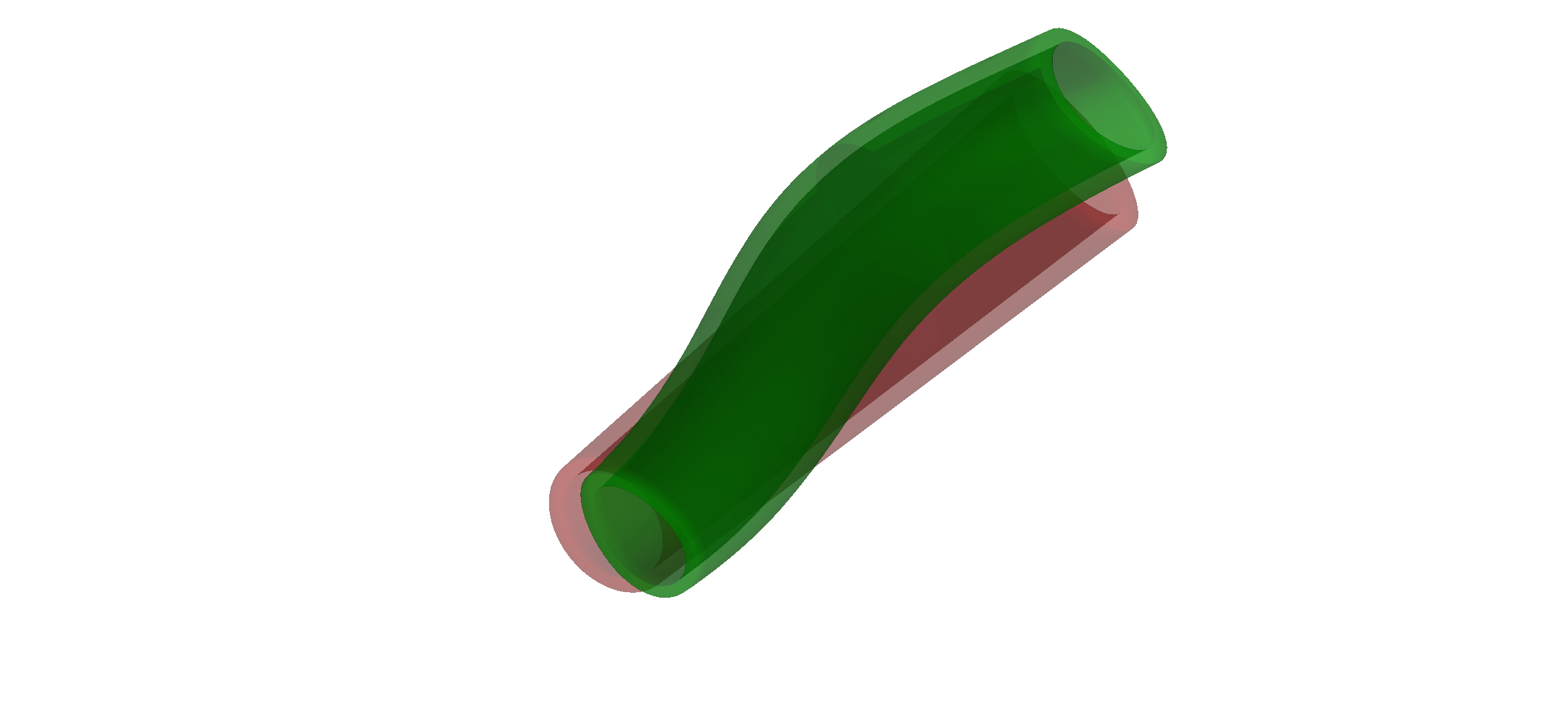}
\includegraphics[width=0.3\textwidth,trim=700pt 50pt 550pt 0pt,clip=true]{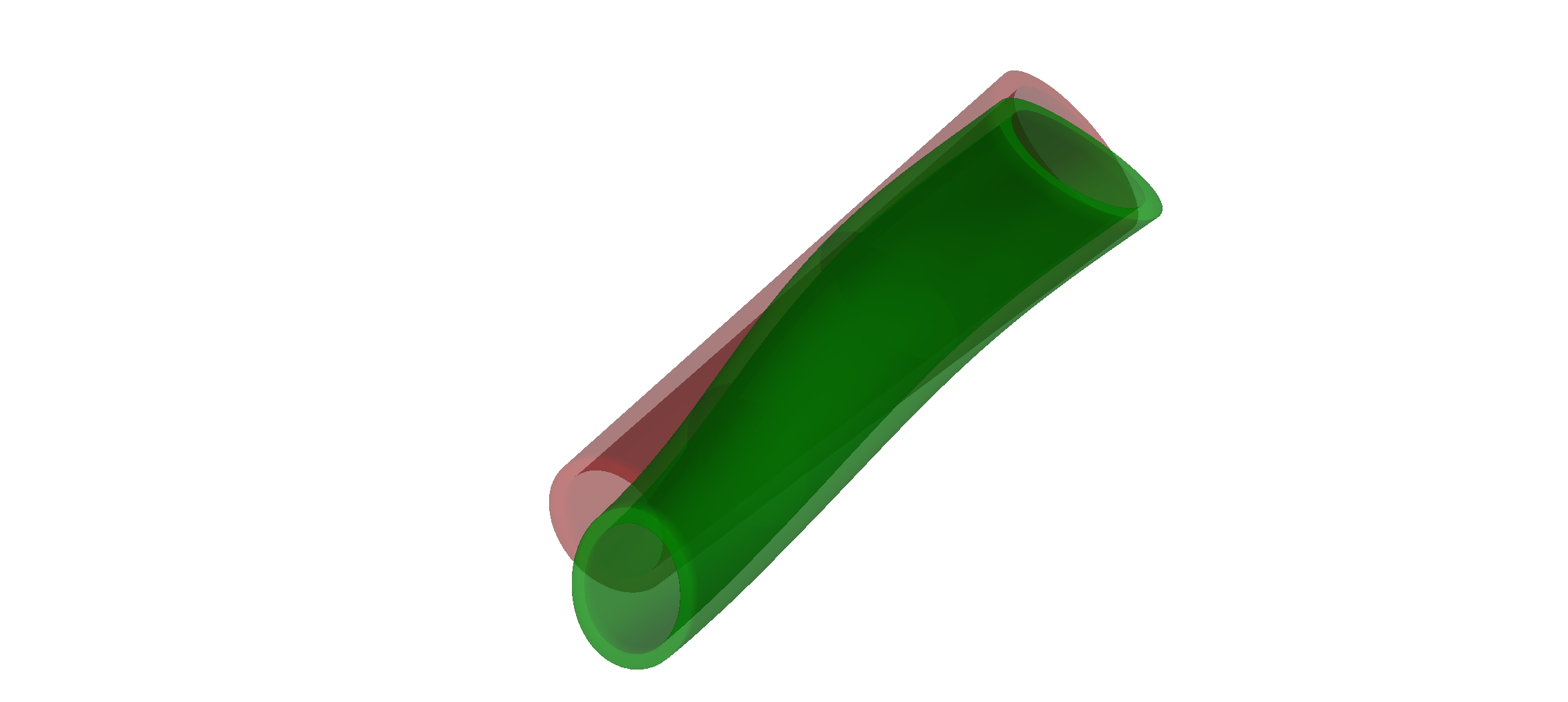}
\includegraphics[width=0.3\textwidth,trim=700pt 50pt 550pt 0pt,clip=true]{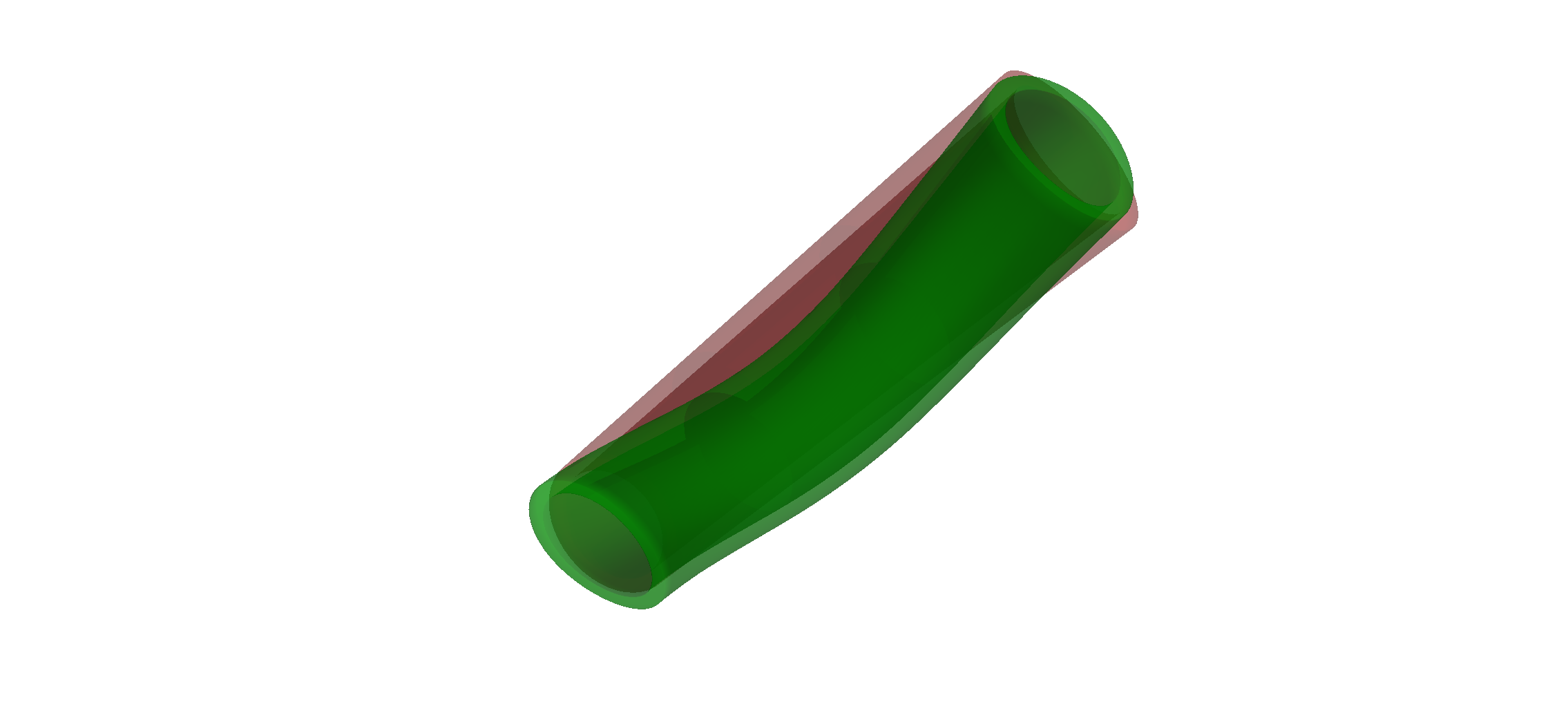}
\includegraphics[width=0.3\textwidth,trim=700pt 100pt 700pt 100pt,clip=true]{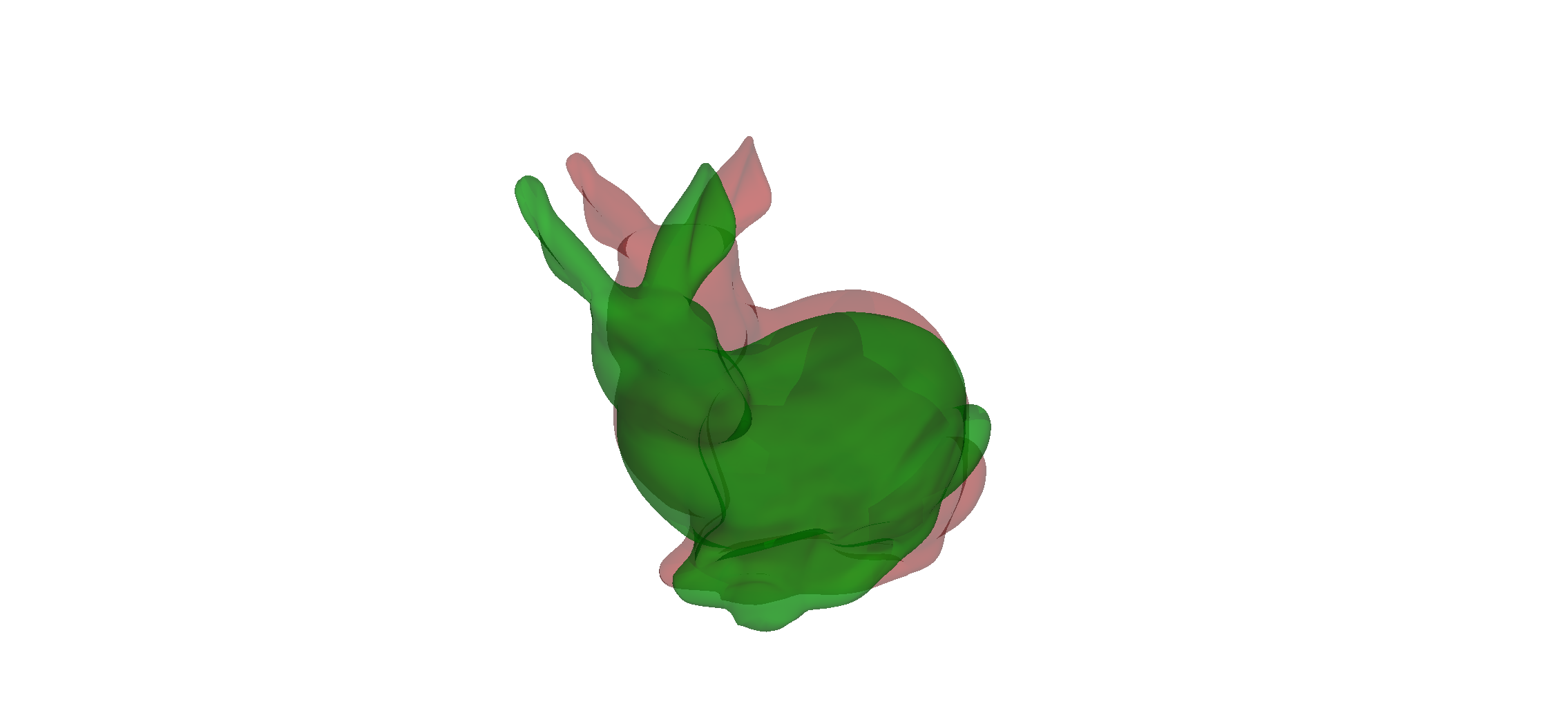}
\includegraphics[width=0.3\textwidth,trim=700pt 100pt 700pt 100pt,clip=true]{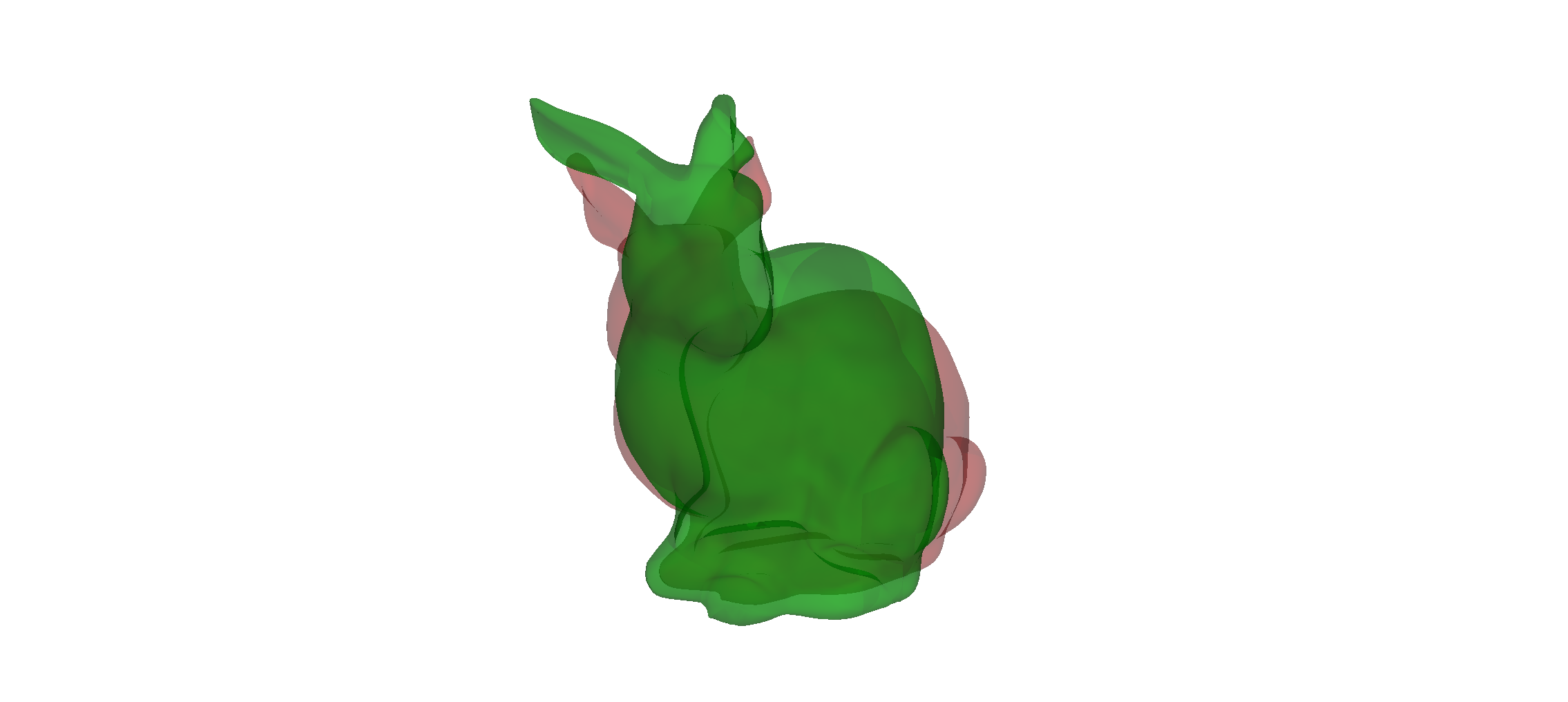}
\includegraphics[width=0.3\textwidth,trim=700pt 100pt 700pt 100pt,clip=true]{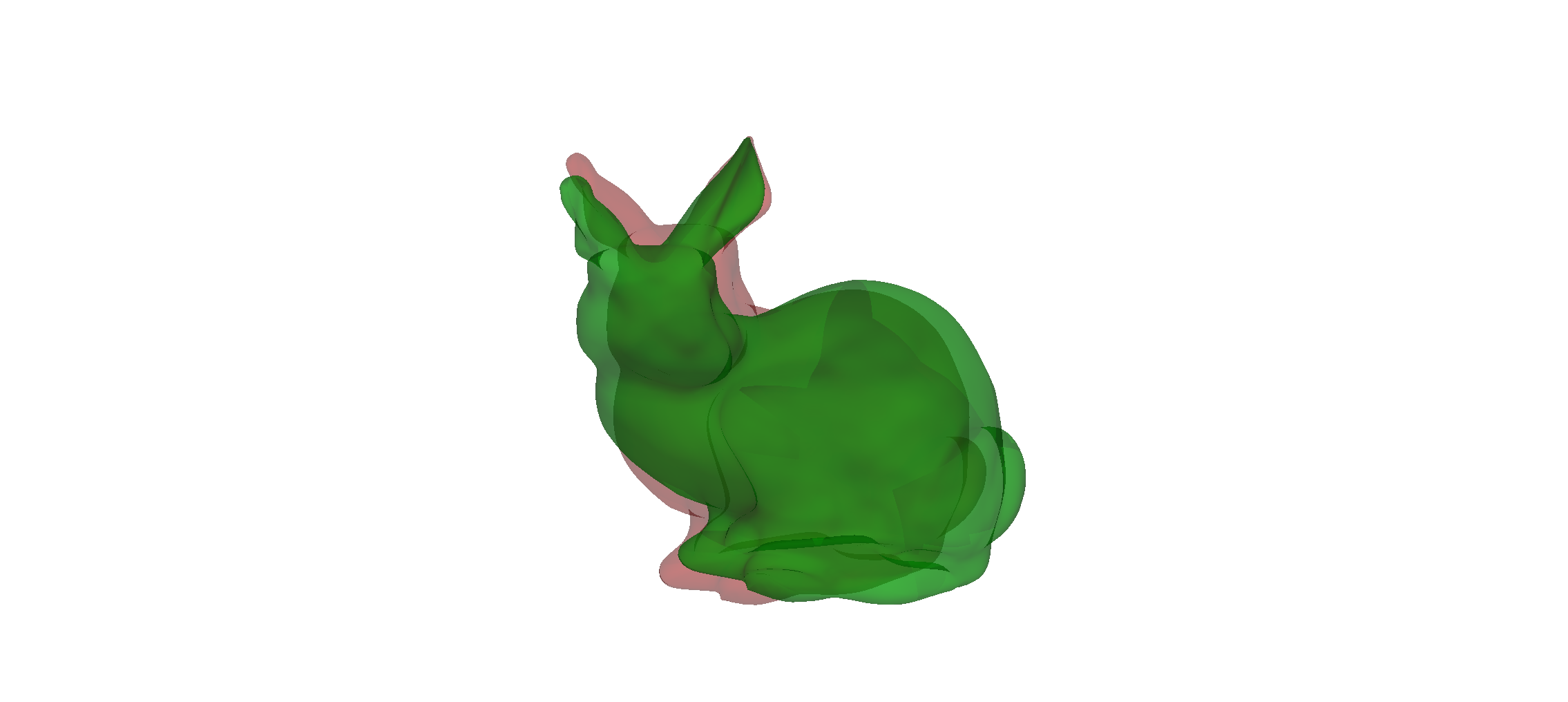}
\caption{\label{fig:def_geo}Different realizations of the random deformation field for the three surfaces
under consideration. The transparent red ones are the reference surfaces, while
green ones are the deformed surfaces.}
\end{center}
\end{figure}

For the random deformation field, we consider the covariance function
\[
\Ebb[{\bs\chi}](\hat{\bs x})=\hat{\bs x},\quad
\Cov[\DefField](\hat{\bs x},\hat{\bs x}')
 =
 10^{-2}
  \begin{bmatrix}
   e^{-50r^2} & 0 & 10^{-4} e^{-0.5r^2} \\
   0 & e^{-50r^2} & 0 \\
   10^{-4}e^{-0.5r^2} & 0 & e^{-50r^2} \\ 
   \end{bmatrix}
\]
with \(r\isdef\|\hat{\bs x}-\hat{\bs x}'\|_2\). The random field is
computed by the pivoted Cholesky decomposition, cf.\ \cite{HPS12,HPS14a,HMvR21},
with an accuracy of \(10^{-4}\). This leads to the parameter dimensions
282, 252, 264 for the sphere, pipe and Stanford bunny, respectively.
Different realizations of the random surfaces
are depicted in Figure~\ref{fig:def_geo}. 

\begin{figure}[htb]
\begin{center}
\includegraphics[scale=0.14,clip,trim=250 20 0 200]{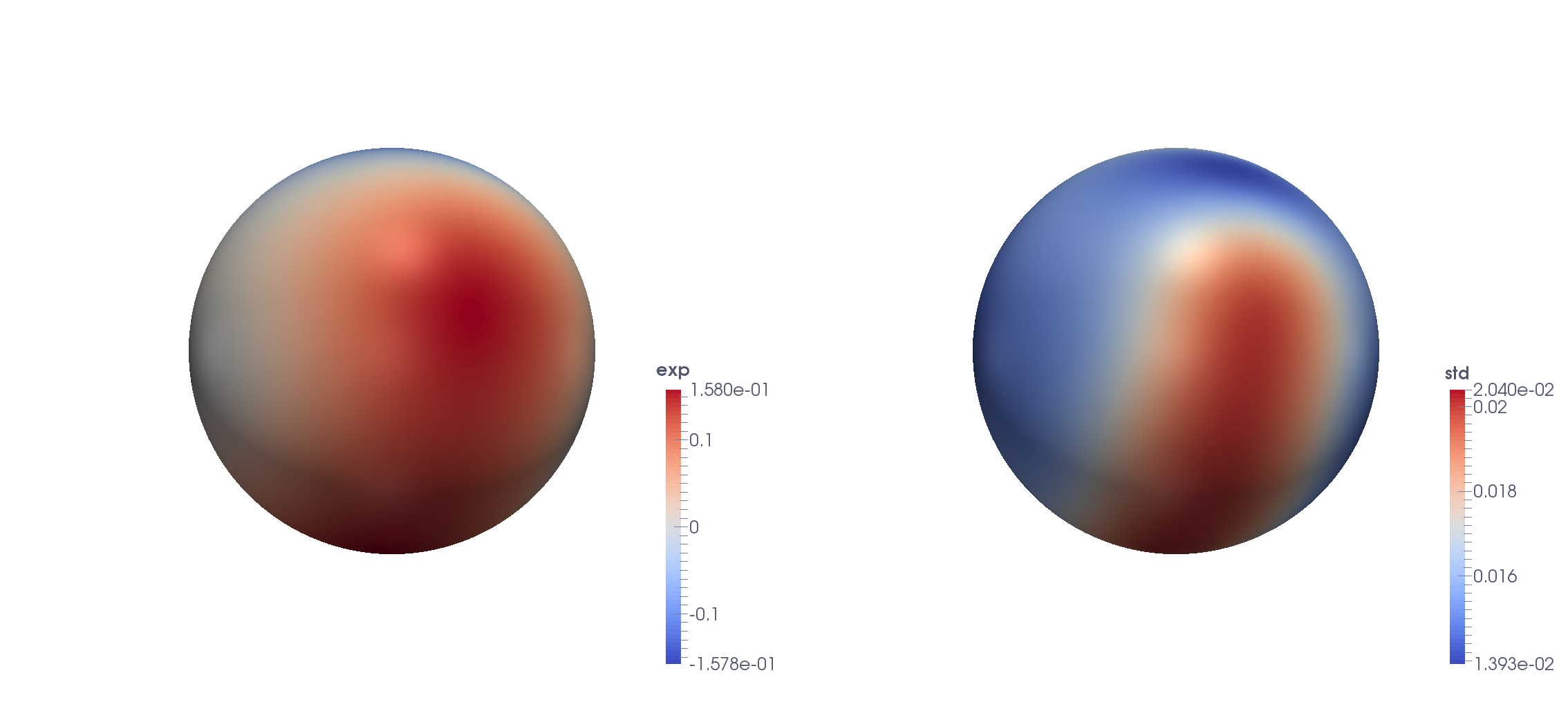}
\includegraphics[width=0.788\textwidth,clip,trim=220 20 0 40]{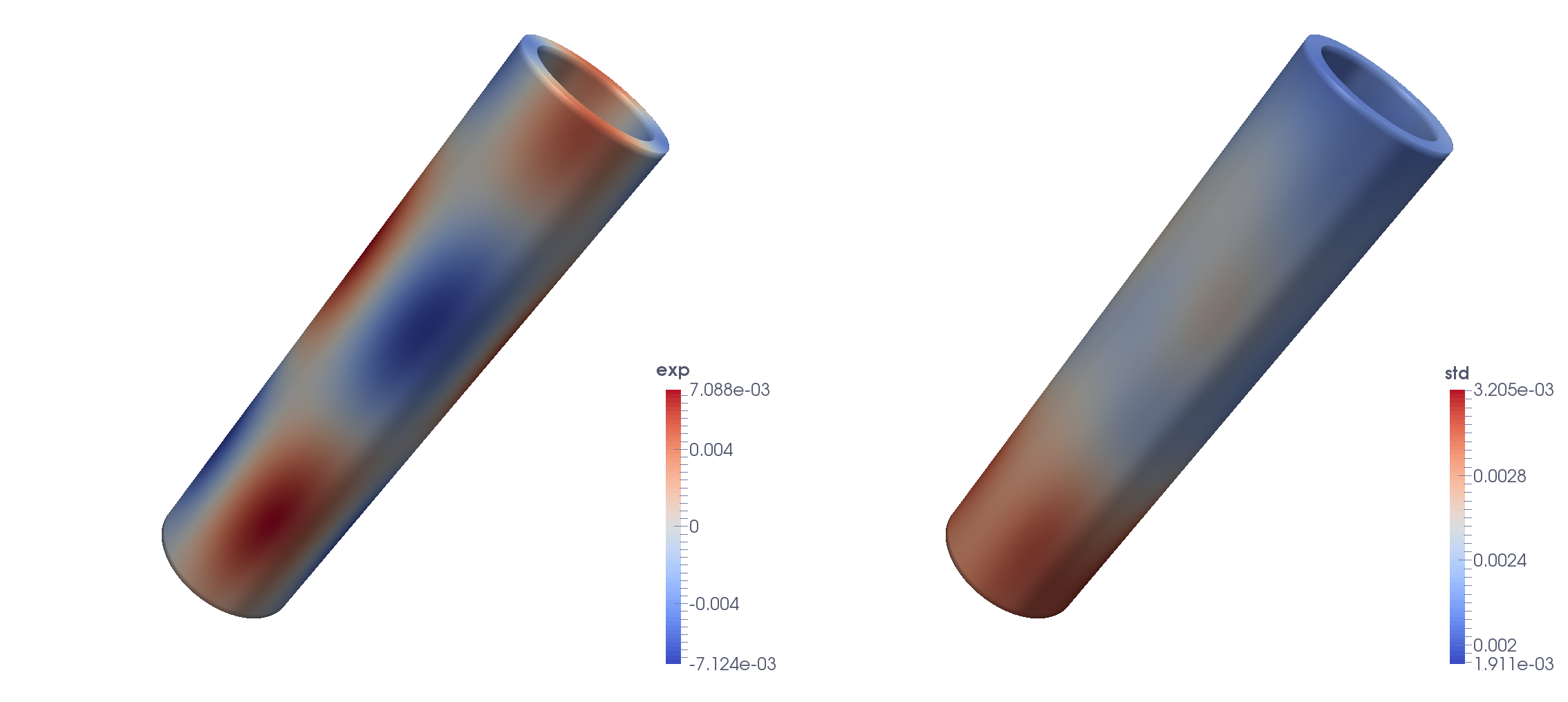}
\includegraphics[width=0.788\textwidth,clip,trim=250 20 0 180]{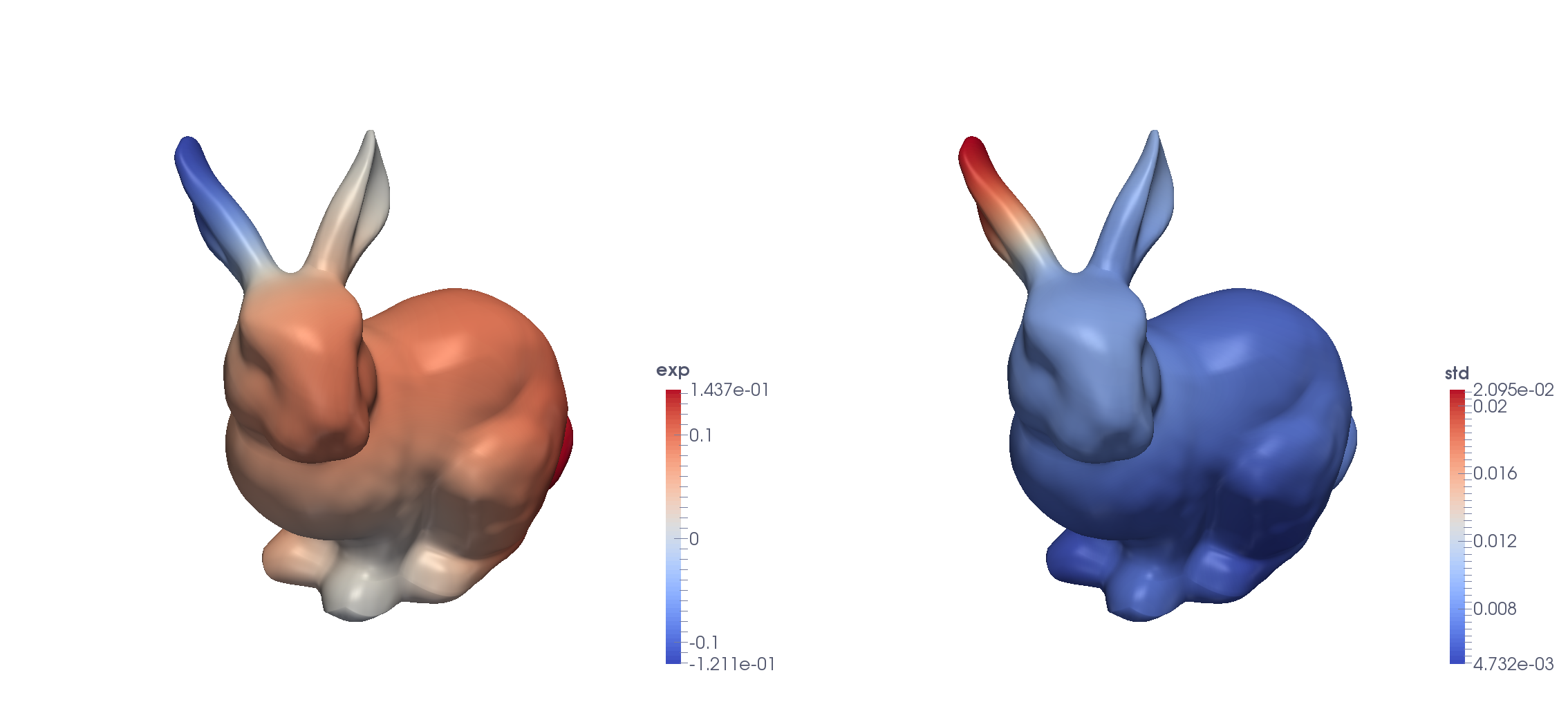}
\caption{\label{fig:bunny_exp_dev}Expectations of the solutions (left)
and standard deviations (right) at 
time \(T=1\) for the three surfaces.}
\end{center}
\end{figure}

The expectations and the standard deviations at the end of time \(T=1\) are
displayed on the reference surfaces in Figure~\ref{fig:bunny_exp_dev}. 

As can be seen, the ratio between the standard deviation and
the expectation is around \(13\%\) for the sphere and the Stanford bunny,
while it is about \(40\%\) for the pipe.

Next, we examine the convergence of the low rank representation. As
we cannot store the full coefficient matrix of the
space-time correlation as a reference,
we only compute the diagonal blocks \({\bs C}_{N_q}(t_i,t_{i})\) for
\(t_i=0.1i\), \(i=1,\ldots,10\), and compare them to the low rank approximation
\(\tilde{{\bs C}}_{N_q}(t_i,t_{i})\). We set the number of samples to
\(N_q = 8192\) for the Monte Carlo and the quasi-Monte Carlo method
based on Halton points. For the assessment of the error, we measure the relative
Frobenius norm of the difference between the ground truth and the low rank
via 
\[
e_{F}\isdef \frac{1}{10}\sum_{i=1}^{10}\frac{\|{\bs C}({t_i},{t_i}) -
\tilde{\bs C}({t_i},{t_i})\|_F}{\|{\bs C}({t_i},{t_i})\|_F}.
\]

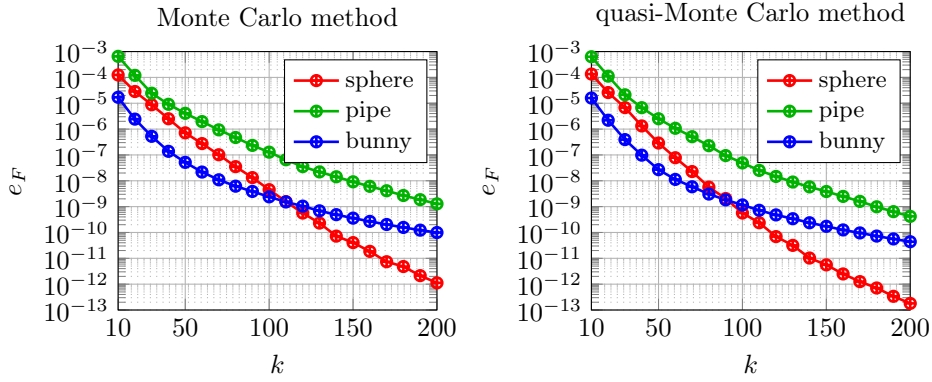
\begin{figure}[htb]
\begin{center}
\pgfplotsset{width=0.46\textwidth}
\pgfplotsset{minor grid style={dotted,black!30}}
\pgfplotsset{grid style={solid,black!30}}
\begin{tikzpicture}
\begin{semilogyaxis}[grid=both,
ymin=1e-13,ymax=1e-3, ytick distance=10, xmin=10,xmax=200,
xtick={10,50,100,150,200},legend style={legend pos=north east,font=\small},
legend cell align={left},%
ylabel={$e_F$}, xlabel ={\(k\)}, title = Monte Carlo method]
\addplot[line width=1pt,color=red,mark=oplus] table [x index = {0}, y index = {1}, col sep=comma] {./Data/sphere/mc_convergence.csv};\addlegendentry{sphere};
\addplot[line width=1pt,color=green!70!black,mark=oplus] table [x index = {0}, y index = {1}, col sep=comma] {./Data/pipe/mc_convergence.csv};\addlegendentry{pipe};
\addplot[line width=1pt,color=blue,mark=oplus] table [x index = {0}, y index = {1}, col sep=comma] {./Data/bunny/mc_convergence.csv};\addlegendentry{bunny};
\end{semilogyaxis}
\end{tikzpicture}
\begin{tikzpicture}
\begin{semilogyaxis}[grid=both, legend style={legend pos=north east,font=\small},
legend cell align={left},%
		ylabel={$e_F$}, xlabel ={\(k\)}, ymin=1e-13,ymax=1e-3, xmin=10,xmax=200,
		ytick distance=10,xtick={10,50,100,150,200},
		title = quasi-Monte Carlo method]
\addplot[line width=1pt,color=red,mark=oplus] table [x index = {0}, y index = {1}, col sep=comma] {./Data/sphere/halton_convergence.csv};\addlegendentry{sphere};
\addplot[line width=1pt,color=green!70!black,mark=oplus] table [x index = {0}, y index = {1}, col sep=comma] {./Data/pipe/halton_convergence.csv};\addlegendentry{pipe};
\addplot[line width=1pt,color=blue,mark=oplus] table [x index = {0}, y index = {1}, col sep=comma] {./Data/bunny/halton_convergence.csv};\addlegendentry{bunny};
\end{semilogyaxis}
\end{tikzpicture}
\caption{\label{fig:convg}Convergence of the low rank approximation for the
		Monte Carlo method and the quasi-Monte Carlo method.}
\end{center}
\end{figure}

The left hand side of Figure~\ref{fig:convg} shows the error \(e_{F}\)
in case of the Monte Carlo method, while the right hand side
provides it for the quasi-Monte Carlo method. In both cases, the error
rapidly decays for increasing \(k\). As can be seen, the errors
for both methods are comparable, however, they are slightly smaller
for the quasi-Monte Carlo method.

\begin{figure}[htb]
\begin{center}
			\includegraphics[clip=true,trim=495 100 440 100, scale=0.14]{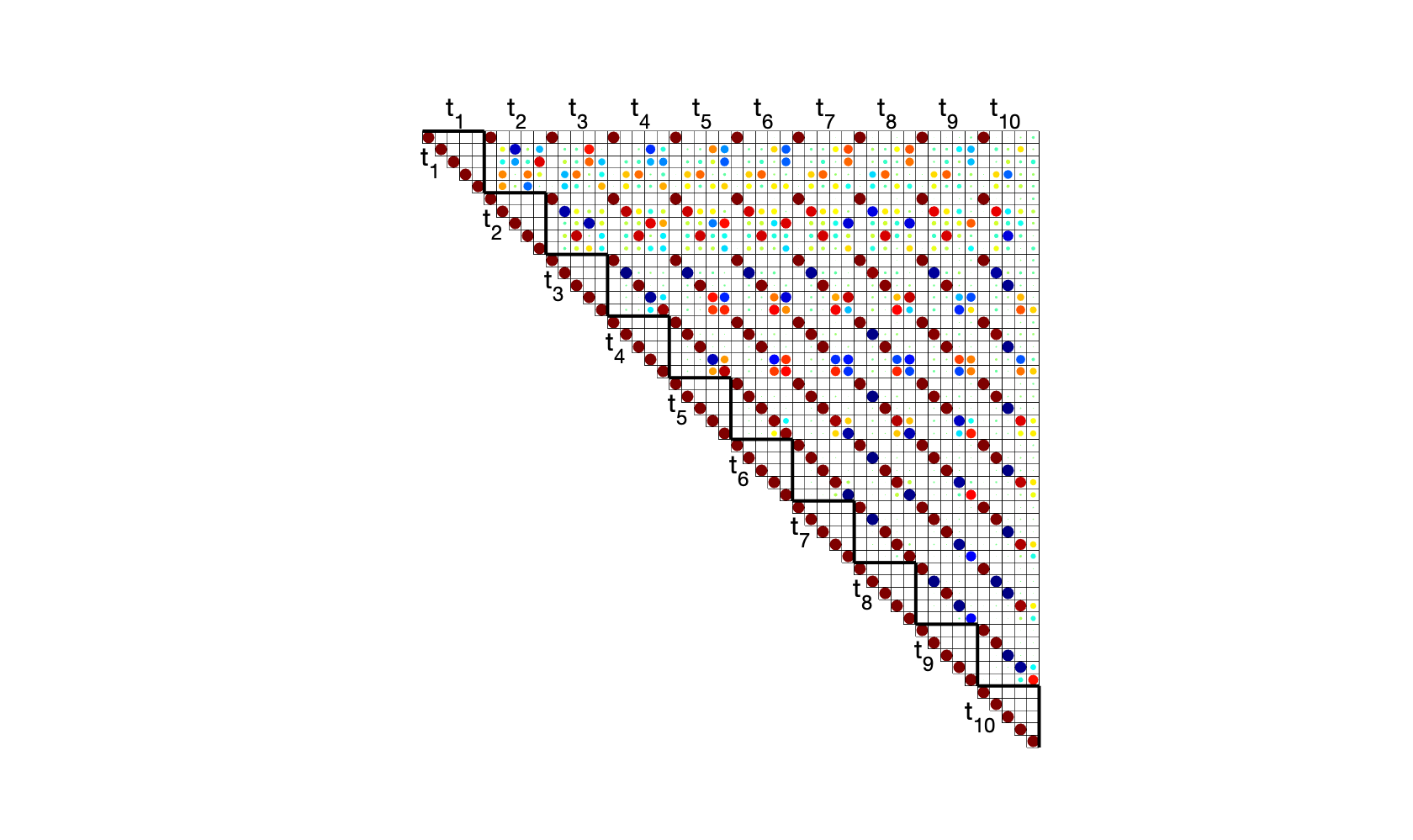}
			\includegraphics[clip=true,trim=495 100 440 100, scale=0.14]{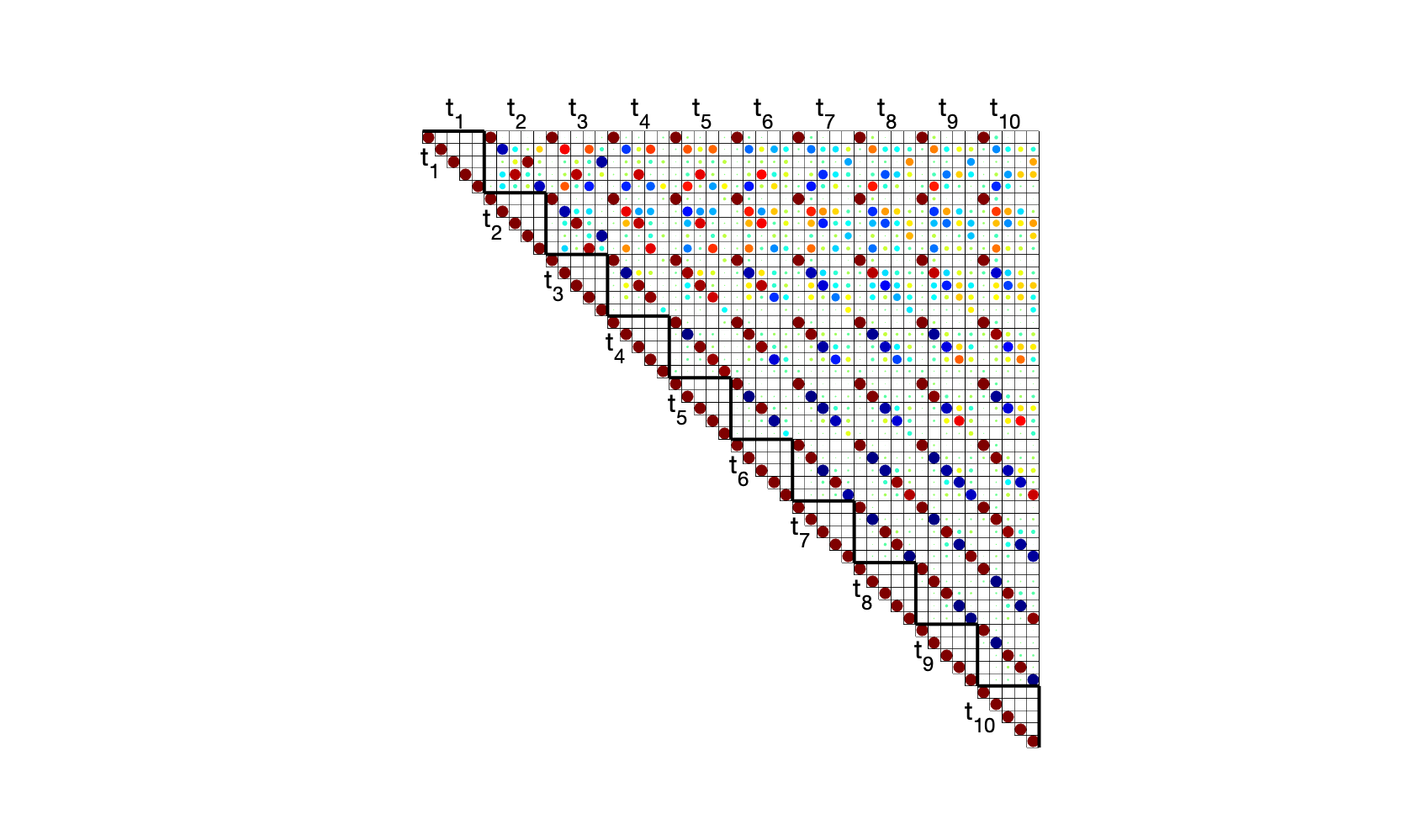}
			\includegraphics[clip=true,trim=415 100 330 100, scale=0.14]{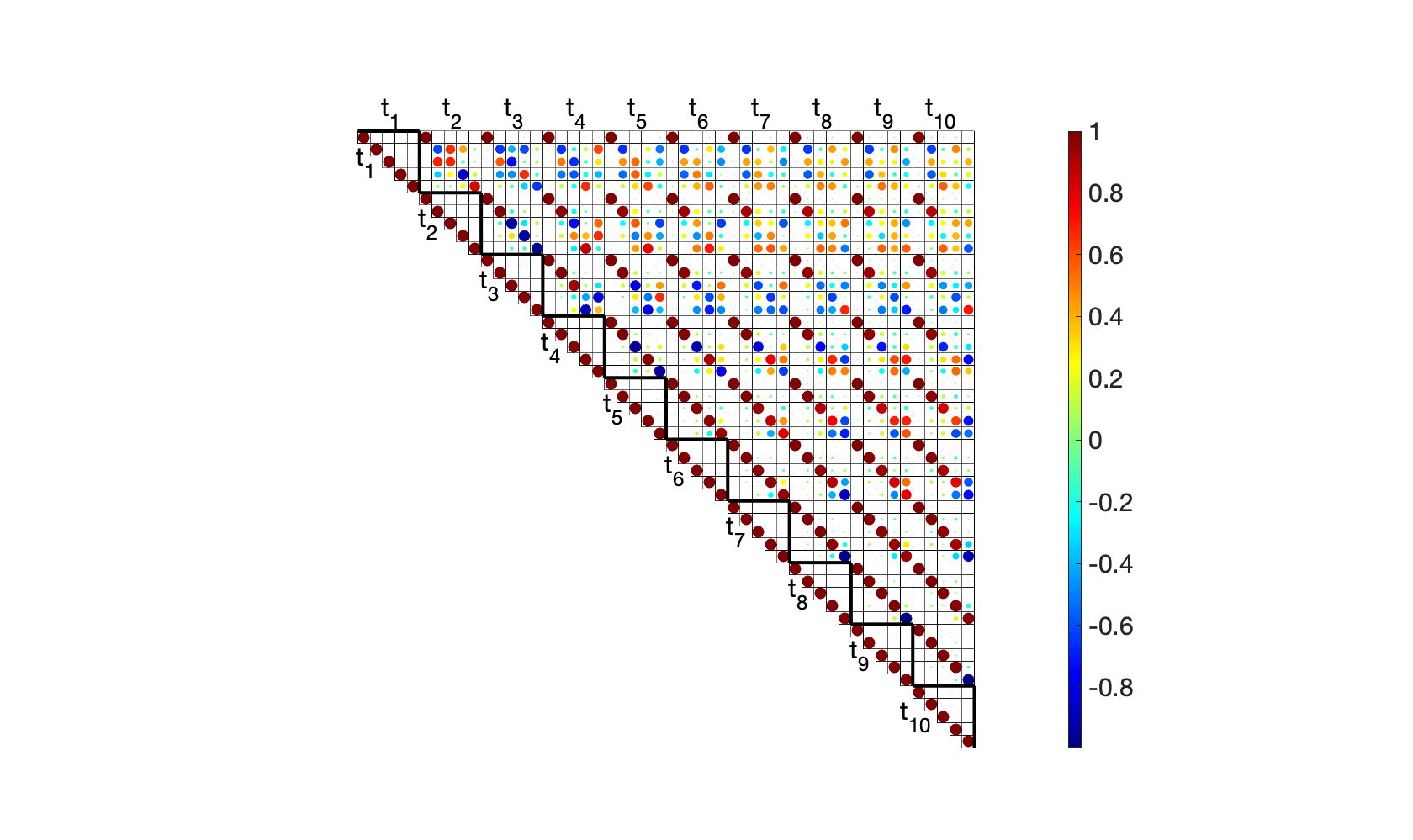}
	\caption{Correlation coefficients for the first 5 principal components
	and the time points \(t_i\) for the sphere (left), the pipe (middle) and
	the Stanford bunny (right) in case of
	the quasi-Monte Carlo method.}
	\label{fig:space-time_cor_qmc}
	\end{center}
\end{figure}

The availability of the full space-time correlation facilitates the visualization
and analysis of the correlation between different time points.
In particular, due to the particular structure of the low rank approximation,
we can directly identify the principal components 
between any two evaluated times, i.e., 
\({\bs U}_\ell^\intercal(t_i)\tilde{{\bs C}}(t_i,t_{i'}){\bs U}_\ell(t_{i'})\),
\(\ell\leq k\). For simplicity, we only show the upper triangular part
of the correlation matrix for the coefficients of the first 5 principal components 
at each time in Figure~\ref{fig:space-time_cor_qmc} for the case of the
quasi-Monte Carlo method. The radii of the circles inside each small block stand for the absolute values of the corresponding correlation coefficients. From the figures, the correlation matrix within the same time step is a \(5 \times 5\) diagonal matrix
(outlined by the black lines in the figures), since we use an orthogonal
basis given by \({\bs U}(t_i)\) at each evaluated time. On the other hand, 
the solution at different times are highly correlated to each other.

\begin{figure}[htb]
\begin{center}
	\pgfplotsset{width=0.46\textwidth}
	\pgfplotsset{minor grid style={dotted,black!30}}
	\pgfplotsset{grid style={solid,black!30}}
		\begin{tikzpicture}
		\begin{semilogyaxis}[grid=both, xmin=10,xmax=200,
		legend style={legend pos=north east,font=\small}, legend cell align={left},%
		ylabel={$e_{\ell^2}$}, xtick={10,50,100,150,200},ymin=1e-6, ymax=1e-1, ytick distance=10,
		xlabel ={\(k\)}, title = 
		Monte Carlo method]
		\addplot[line width=1pt,color=red,mark=oplus] table [x index = {0}, y index = {1}, col sep=comma] {./Data/sphere/mc_generalization.csv};\addlegendentry{sphere};
		\addplot[line width=1pt,color=green!70!black,mark=oplus] table [x index = {0}, y index = {1}, col sep=comma] {./Data/pipe/mc_generalization.csv};\addlegendentry{pipe};
		\addplot[line width=1pt,color=blue,mark=oplus] table [x index = {0}, y index = {1}, col sep=comma] {./Data/bunny/mc_generalization.csv};\addlegendentry{bunny};
		\end{semilogyaxis}
		\end{tikzpicture}
		\begin{tikzpicture}
		\begin{semilogyaxis}[xmin=10,xmax=200,
		grid=both, legend style={legend pos=north east,font=\small}, legend cell align={left},%
		ylabel={$e_{\ell^2}$},ytick distance=10, ymin=1e-6, ymax=1e-1,
		xtick={10,50,100,150,200},xlabel ={\(k\)}, title = quasi-Monte Carlo method]
		\addplot[line width=1pt,color=red,mark=oplus] table [x index = {0}, y index = {1}, col sep=comma] {./Data/sphere/halton_generalization.csv};\addlegendentry{sphere};
		\addplot[line width=1pt,color=green!70!black,mark=oplus] table [x index = {0}, y index = {1}, col sep=comma] {./Data/pipe/halton_generalization.csv};\addlegendentry{pipe};
		\addplot[line width=1pt,color=blue,mark=oplus] table [x index = {0}, y index = {1}, col sep=comma] {./Data/bunny/halton_generalization.csv};\addlegendentry{bunny};
		\end{semilogyaxis}
		\end{tikzpicture}
\caption{\label{fig:gen}Generalization experiments on 8192 test samples generated by the Monte Carlo method.}
\end{center}
\end{figure}
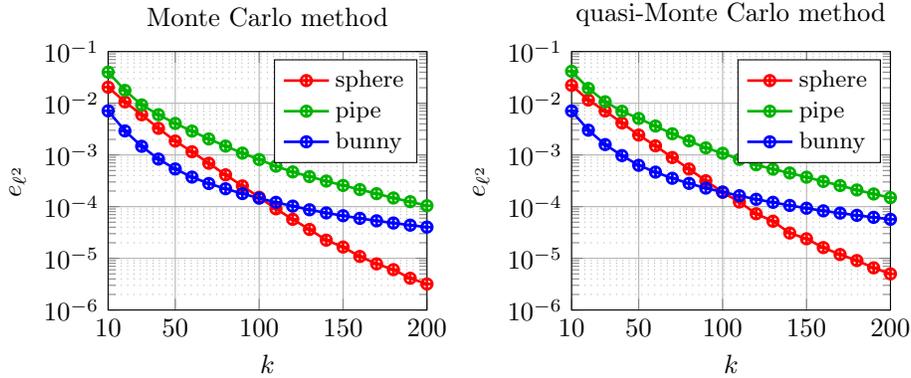

Next, we consider the generalization error given the low rank bases
 \({\bs U}(t_i)\).
 To this end, we introduce the relative
 \(\ell^2\)-error at the time point \(t_i\), i.e.,
\[
e_{\ell^2}(t_i) = \frac{1}{N_q} \sum_{j=1}^{N_q}
\frac{\big\|{\bs u}(\tilde{\bs\xi}_j,t_i) - {\bs U}(t_i){\bs U}(t_i)^\intercal{\bs u}(\tilde{\bs\xi}_j,t_i)\big\|_2}{\big\|{\bs u}(\tilde{\bs\xi}_j,t_i)\big\|_2},
\]
for independently drawn Monte Carlo
samples \({\bs u}(\tilde{\bs\xi}_j,t_i)\), \(j=1,\ldots N_q\),
with \(N_q=8192\). As a measure for the generalization error, we
consider now the average error
 \[
e_{\ell^2} = \frac{1}{N_t} \sum_{i=1}^{N_t} e_{\ell^2}(t_i).
\]
The generalization error is displayed in Figure~\ref{fig:gen}.
Particularly for the sphere, the generalization error rapidly decreases
for an increasing rank \(k\), where the Monte Carlo method performs slightly
better here. For \(k=200\) the average \(\ell^2\)-error 
is lowest for the sphere and around \(4\cdot 10^{-6}\) for both methods.

Finally, we examine the stability of the low rank approximation with respect
to the magnitude of the random perturbation. To this end, we introduce
the parameter \(\alpha\) which steers the impact of the random perturbation
according to
\begin{equation*}
\DefField({\bs y},\xref)=\Ebb[\DefField](\xref)+
\alpha \sum_{k=1}^m\sqrt{\lambda_{k}}\DefField_{k}(\xref)y_k,
\quad{\bs y}\in\Gamma\isdef[-1,1]^{m}.
\end{equation*}

\begin{figure}[htb]
\begin{center}
\pgfplotsset{width=0.46\textwidth}
\pgfplotsset{minor grid style={dotted,black!30}}
\pgfplotsset{grid style={solid,black!30}}
\begin{tikzpicture}
\begin{semilogyaxis}[
ytick distance=10,xtick distance=0.2,xmin=1,xmax=2,ymin=1e-14,ymax=1e-8,
grid=both, legend style={legend pos=south east,font=\small, title = Monte Carlo method}, legend cell align={left},%
			ylabel={$e_F$}, xlabel ={$\alpha$}]
			\addplot[line width=1pt,color=red,mark=oplus] table [x index = {0}, y index = {1}, col sep=comma] {./Data/sphere/mc_stability_cov.csv};\addlegendentry{sphere};
			\addplot[line width=1pt,color=green!70!black,mark=oplus] table [x index = {0}, y index = {1}, col sep=comma] {./Data/pipe/mc_stability_cov.csv};\addlegendentry{pipe};
			\addplot[line width=1pt,color=blue,mark=oplus] table [x index = {0}, y index = {1}, col sep=comma] {./Data/bunny/mc_stability_cov.csv};\addlegendentry{bunny};
			\end{semilogyaxis}
			\end{tikzpicture}
			\begin{tikzpicture}
			\begin{semilogyaxis}[
			ytick distance=10,xtick distance=0.2,xmin=1,xmax=2,ymin=1e-14,ymax=1e-8,
			grid=both, legend style={legend pos=south east,font=\small}, legend cell align={left},%
			ylabel={$e_F$}, xlabel ={$\alpha$},title = quasi-Monte Carlo method]
			\addplot[line width=1pt,color=red,mark=oplus] table [x index = {0}, y index = {1}, col sep=comma] {./Data/sphere/halton_stability_cov.csv};\addlegendentry{sphere};
			\addplot[line width=1pt,color=green!70!black,mark=oplus] table [x index = {0}, y index = {1}, col sep=comma] {./Data/pipe/halton_stability_cov.csv};\addlegendentry{pipe};
			\addplot[line width=1pt,color=blue,mark=oplus] table [x index = {0}, y index = {1}, col sep=comma] {./Data/bunny/halton_stability_cov.csv};\addlegendentry{bunny};
			\end{semilogyaxis}
			\end{tikzpicture}
	\caption{\label{fig:stab_cov}Low rank approximation error with respect to the 
	magnitude of the random perturbation.}
		\end{center}
\end{figure}
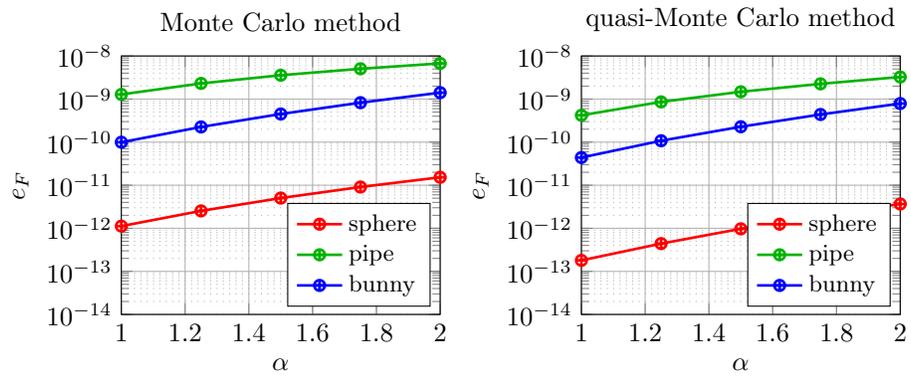

We fix \(k=200\) and consider \(\alpha=1,1.2,\ldots,2\). 
This results in a maximal possible relative displacement between
\(20\%\) and \(40\%\) for the sphere, a maximal possible relative displacement between \(21\%\) and \(42\%\) for the pipe, and a maximal possible relative displacement between \(27\%\) and \(54\%\) for the Stanford bunny.
Figure~\ref{fig:stab_cov}
depicts the corresponding error \(e_F\) for the covariance approximation.

As can be seen, both methods are rather robust with respect to the increase
of the displacement, where the quasi-Monte Carlo method behaves slightly better.

\begin{figure}[htb]
\begin{center}
\pgfplotsset{width=0.46\textwidth}
\pgfplotsset{minor grid style={dotted,black!30}}
\pgfplotsset{grid style={solid,black!30}}
\begin{tikzpicture}
\begin{semilogyaxis}[
ytick distance=10,xtick distance=0.2,xmin=1,xmax=2,ymin=1e-7,ymax=1e-3,
grid=both, legend style={legend pos=south east,font=\small, title = Monte Carlo method}, legend cell align={left},%
			ylabel={$e_{\ell^2}$}, xlabel ={$\alpha$}]
			\addplot[line width=1pt,color=red,mark=oplus] table [x index = {0}, y index = {1}, col sep=comma] {./Data/sphere/mc_stability_gen.csv};\addlegendentry{sphere};
			\addplot[line width=1pt,color=green!70!black,mark=oplus] table [x index = {0}, y index = {1}, col sep=comma] {./Data/pipe/mc_stability_gen.csv};\addlegendentry{pipe};
			\addplot[line width=1pt,color=blue,mark=oplus] table [x index = {0}, y index = {1}, col sep=comma] {./Data/bunny/mc_stability_gen.csv};\addlegendentry{bunny};
			\end{semilogyaxis}
			\end{tikzpicture}
			\begin{tikzpicture}
			\begin{semilogyaxis}[
			ytick distance=10,xtick distance=0.2,xmin=1,xmax=2,ymin=1e-7,ymax=1e-3,
			grid=both, legend style={legend pos=south east,font=\small}, legend cell align={left},%
			ylabel={$e_{\ell^2}$}, xlabel ={$\alpha$},title = quasi-Monte Carlo method]
			\addplot[line width=1pt,color=red,mark=oplus] table [x index = {0}, y index = {1}, col sep=comma] {./Data/sphere/halton_stability_gen.csv};\addlegendentry{sphere};
			\addplot[line width=1pt,color=green!70!black,mark=oplus] table [x index = {0}, y index = {1}, col sep=comma] {./Data/pipe/halton_stability_gen.csv};\addlegendentry{pipe};
			\addplot[line width=1pt,color=blue,mark=oplus] table [x index = {0}, y index = {1}, col sep=comma] {./Data/bunny/halton_stability_gen.csv};\addlegendentry{bunny};
			\end{semilogyaxis}
			\end{tikzpicture}
	\caption{\label{fig:stab_gen}Generalization error with respect to the magnitude of the randomness.}
		\end{center}
\end{figure}
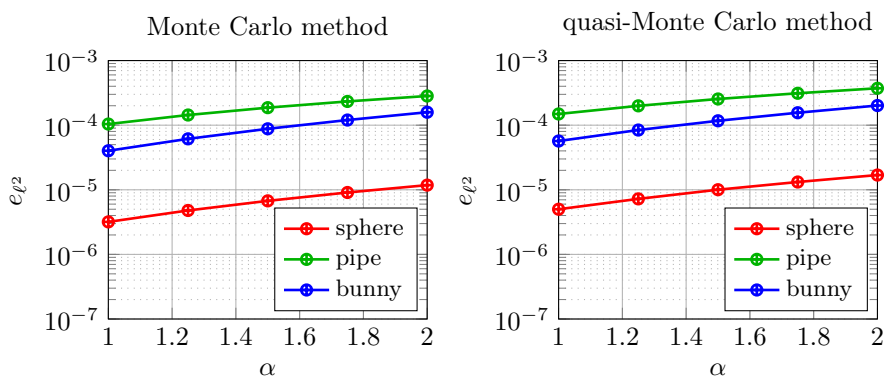

The generalization error with respect to \(\alpha\) is shown in 
Figure~\ref{fig:stab_gen}.

As for the low rank error, the generalization error only moderately increases with
an increase of the random perturbation, where here the Monte Carlo method performs
slightly better.

\section{Conclusion}\label{sec:Conclusion}
In this article, we have presented an isogeometric approach
to solving diffusion problems on random surfaces.
Especially, we have described in detail,
how diffusion problems on random surfaces
can be modelled by means of random deformation fields
in the isogeometric context and how quantities of interest 
may be derived. 
Moreover, we have employed an online low rank approximation
algorithm for the high-dimensional
space-time correlation of the random solution.
Extensive numerical studies on complex geometries have been
performed. The numerical
results corroborate the efficacy of the presented methodology.
Finally, the solver for the Laplace-Beltrami operator will be
added to the isogeometric boundary element library
\verb+Bembel+. 

\bibliography{literature}
\bibliographystyle{abbrv}
\end{document}